\documentclass[a4paper]{amsart}
\usepackage[utf8]{inputenc}
\usepackage[T1]{fontenc}
\usepackage{lmodern}
\usepackage{amssymb}
\usepackage[all]{xy}
\usepackage{nicefrac,mathtools,enumitem}
\usepackage{microtype}

\usepackage[pdftitle={Square-integrable coactions of quantum group on Hilbert modules},
  pdfauthor={Alcides Buss and Ralf Meyer},
  pdfsubject={Mathematics; MSC 19K35, 46L80}
]{hyperref}
\usepackage[lite]{amsrefs}
\newcommand*{\MRref}[2]{ \href{http://www.ams.org/mathscinet-getitem?mr=#1}{MR #1}}
\newcommand*{\arxiv}[1]{ \href{http://www.arxiv.org/abs/#1}{arXiv:#1}}

\numberwithin{equation}{section}
\theoremstyle{plain}
\newtheorem{theorem}[equation]{Theorem}
\newtheorem{lemma}[equation]{Lemma}
\newtheorem{proposition}[equation]{Proposition}
\newtheorem{corollary}[equation]{Corollary}
\theoremstyle{definition}
\newtheorem{definition}[equation]{Definition}

\theoremstyle{remark}
\newtheorem{remark}[equation]{Remark}
\newtheorem{example}[equation]{Example}

\DeclareMathOperator*{\dom}{dom}
\DeclareMathOperator*{\Mor}{Mor}
\DeclareMathOperator*{\Sub}{Sub}
\DeclareMathOperator*{\spann}{span}
\DeclareMathOperator*{\clspan}{\overline{span}}

\newcommand*{\nb}{\nobreakdash}
\newcommand*{\Star}{\texorpdfstring{$^*$\nobreakdash-\hspace{0pt}}{*-}}

\newcommand*{\C}{\mathbb C}

\newcommand*{\R}{\mathbb R}
\newcommand*{\N}{\mathbb N}

\newcommand*{\Bound}{\mathbb B}
\newcommand*{\Comp}{\mathbb K}

\newcommand*{\dd}{\textup d}
\newcommand*{\ima}{\textup i}
\newcommand*{\si}{\textup{si}}
\newcommand*{\ible}{\textup{i}}
\newcommand*{\comm}{\textup{c}}
\newcommand*{\red}{\textup{r}}
\newcommand*{\uni}{\textup{u}}
\newcommand*{\opp}{\textup{op}}

\newcommand*{\CONT}{\mathcal C}
\newcommand*{\ID}{\textup{id}}

\newcommand*{\Cred}{\textup{C}^*_\textup{r}}
\newcommand*{\QG}{\mathcal{G}}
\newcommand*{\com}{\Delta}
\newcommand*{\Haar}{\varphi}
\newcommand*{\coact}{\Delta}
\newcommand*{\Mult}{\mathcal M}
\newcommand*{\dual}{\ast}
\newcommand*{\Hils}{\mathcal H}
\newcommand*{\Hilsr}{\mathcal L}
\newcommand*{\Hilm}{\mathcal E}
\newcommand*{\modauto}{\sigma}
\newcommand*{\modular}{\delta}
\newcommand*{\QGtodual}{\mu}
\newcommand*{\coef}{B}
\newcommand*{\coefd}{D}
\newcommand*{\Av}{\textup{Av}}

\newcommand*{\defeq}{\mathrel{\vcentcolon=}}
\newcommand*{\eqdef}{\mathrel{=\vcentcolon}}
\newcommand*{\into}{\rightarrowtail}
\newcommand*{\prto}{\twoheadrightarrow}

\newcommand*{\abs}[1]{\lvert#1\rvert}
\newcommand*{\norm}[1]{\lVert#1\rVert}
\newcommand*{\conj}[1]{\overline{#1}}
\newcommand{\ket}[1]{\mathopen|#1\rangle}
\newcommand{\bra}[1]{\langle#1\mathclose|}
\newcommand{\KET}[1]{\mathopen|#1\rangle\!\rangle}
\newcommand{\BRA}[1]{\langle\!\langle#1\mathclose|}

\begin{document}
\title[Square-integrable coactions]{Square-integrable coactions\\
  of locally compact quantum groups}

\author{Alcides Buss}
\email{alcides@mtm.ufsc.br}
\address{Departamento de Matemática\\
  Universidade Federal de Santa Catarina\\
  88.040-900 Florianópolis-SC\\
  Brasil}

\author{Ralf Meyer}
\email{rameyer@uni-math.gwdg.de}
\address{Mathematisches Institut\\
  Georg-August-Universität Göttingen\\
  Bunsenstraße 3--5\\
  37073 Göttingen\\
  Germany}

\begin{abstract}
  We define and study square-integrable coactions of locally
  compact quantum groups on Hilbert modules, generalising
  previous work for group actions.  As special cases, we
  consider square-integrable Hilbert space corepresentations
  and integrable coactions on \(C^*\)\nb-algebras.  Our main
  result is an equivariant generalisation of Kasparov's
  Stabilisation Theorem.
\end{abstract}


\thanks{This research was supported by the EU-Network
  \emph{Quantum Spaces and Noncommutative Geometry} (Contract
  HPRN-CT-2002-00280), the Marie Curie Action
  \emph{Noncommutative Geometry and Quantum Groups} (Contract
  MKTD-CT-2004-509794), the \emph{Deutsche
    Forschungsgemeinschaft} (SFB 478), and by CAPES and CNPq/PNPD Process Number 558420/2008-7 (Brazil).
    Part of this work is based on the first author's doctoral thesis
    under the supervision of the second author and Siegfried Echterhoff.}

\maketitle

\section{Introduction}
\label{sec:introduction}

This article generalises previous work for group actions by
Marc Rieffel and the second author in
\cites{Rieffel:Integrable_proper, Meyer:Generalized_Fixed} to
coactions of locally compact quantum groups.  First we briefly
explain the results for group actions we are going to
generalise.

Square-integrable group representations have played an
important role in representation theory for a long time.
Roughly speaking, a unitary representation~\(\pi\) of a locally
compact group~\(G\) on a Hilbert space~\(\Hils\) is
square-integrable if there are enough vectors
\(\xi,\eta\in\Hils\) for which the function \(c_{\xi\eta}(g)
\defeq \langle \pi_g(\xi),\eta\rangle\) belongs to the Hilbert
space \(L^2(G)\), defined using a left invariant Haar measure
on~\(G\).  Irreducible square-integrable representations behave
in many respects like irreducible representations of compact
groups: they are strongly contained in the regular
representation, and the coefficients~\(c_{\xi\eta}\) satisfy
orthogonality relations reminiscent of the relations in the
Peter--Weyl Theorem---except that the modular function
complicates the formulas for non-unimodular groups.

The definition of~\(c_{\xi\eta}\) makes perfect sense if
\(\xi\) and~\(\eta\) are elements of a Hilbert
\(\coef\)\nb-module, where~\(\coef\) is some \(C^*\)\nb-algebra
with a continuous action of~\(G\).  Square-integrability
requires this coefficient to belong to the Hilbert
\(\coef\)\nb-module \(L^2(G,\coef)\).  The main general result
about square-integrable group actions on Hilbert modules is
that a countably generated \(G\)\nb-equivariant Hilbert
\(\coef\)\nb-module~\(\Hilm\) is square-integrable if and only
if
\[
\Hilm\oplus L^2(G,\coef)^\infty \cong L^2(G,\coef)^\infty.
\]
This is an equivariant analogue of Kasparov's Stabilisation
Theorem and generalises the relationship between
square-integrable Hilbert space representations and the regular
representation.

A group action~\(\beta\) on a \(C^*\)\nb-algebra~\(\coef\) is
integrable if there is a dense subspace of positive elements
\(x\in\coef^+\) for which the integral \(\int_G \beta_g(x)\,\dd
g\) converges (unconditionally) in the strict topology.  This
notion is closely related to square-integrability.  A group
action on a Hilbert module~\(\Hilm\) is square-integrable if
and only if the induced action on \(\Comp(\Hilm)\) is
integrable.  As a consequence, an action on a
\(C^*\)\nb-algebra~\(\coef\) is integrable if and only
if~\(\coef\) is square-integrable as an equivariant Hilbert
module over itself.

An action on a commutative \(C^*\)\nb-algebra \(\CONT_0(X)\) is
integrable if and only if the corresponding action on~\(X\) is
proper (see~\cite{Rieffel:Integrable_proper}).  But for
non-commutative \(C^*\)\nb-algebras, integrability is more
closely related to stability than to properness.  The
equivariant stabilisation theorem forces the diagonal action on
\(\Hilm\otimes L^2(G)\) to be square-integrable for any Hilbert
module~\(\Hilm\), so that the diagonal action on
\(\coef\otimes\Comp(L^2G)\) is integrable for any
\(C^*\)\nb-algebra~\(\coef\).  Furthermore, dual actions and
actions of compact groups are integrable.

In this article, we are going to extend these results to
coactions of locally compact quantum groups. Group actions
become coactions of the locally compact quantum group
\((\CONT_0(G),\com)\) with \(\com(f)(x,y)\defeq f(x\cdot y)\)
for all \(x,y\in G\), \(f\in\CONT_0(G)\).

The main new difficulty is that the straightforward
relationship between \(\CONT_0(G)\), \(L^2(G)\), and
\(\Cred(G)\) becomes less transparent.  Let \((\QG,\com)\)
be a locally compact quantum group.  Then~\(\QG\) corresponds
to \(\CONT_0(G)\).  The analogue of \(L^2(G)\) is the Hilbert
space~\(\Hilsr\) associated to a left invariant Haar weight
on~\(\QG\).  It is related to~\(\QG\) by a densely defined,
closed, unbounded linear map with dense range
\[
\Lambda\colon \QG\supseteq \dom(\Lambda) \to \Hilsr,
\]
which is characterised by the condition \(\bigl\langle
\Lambda(x),\Lambda(y)\bigr\rangle = \Haar(x^*y)\) for all
\(x,y\in\dom(\Lambda)\), where~\(\Haar\) denotes the left
invariant Haar weight on~\(\QG\).  In the group case,
\(\Lambda\) is the “identical” map that views a function in
\(\dom(\Lambda) = \CONT_0(G)\cap L^2(G)\subseteq\CONT_0(G)\) as
an element of \(L^2(G)\).  The reduced group \(C^*\)\nb-algebra
\(\Cred(G)\) corresponds to the dual quantum
group~\(\hat{\QG}\).  This has its own left invariant Haar
weight~\(\hat{\Haar}\), and the corresponding Hilbert space may
be identified with~\(\Hilsr\).  This involves another densely
defined, closed, unbounded linear map with dense range
\[
\hat{\Lambda}\colon \hat{\QG}\supseteq \dom(\hat{\Lambda}) \to
\Hilsr.
\]
In the group case, \(\dom(\hat{\Lambda})\) contains
\(L^1(G)\cap L^2(G)\), viewed as a subspace of \(\Cred(G)\),
and maps it “identically” to \(L^2(G)\).

Using the unbounded linear maps \(\Lambda\)
and~\(\hat{\Lambda}\), we carry over the theory of
square-integrable Hilbert space representations of groups to
corepresentations of \((\QG,\com)\) on Hilbert spaces in
Section~\ref{sec:si_Hilbert_space}.  To prove that these
corepresentations have the expected properties, it is
convenient to view them as representations of the (universal)
dual quantum group of \((\QG,\com)\).  The resulting notion of
square-integrable representation is a special case of a
definition by François
Combes~\cite{Combes:Systemes_Hilbertiens}.

Our proof of the orthogonality relations for irreducible
square-integrable corepresentations
follows~\cite{Rieffel:Integrable_proper}: the crucial
ingredient is the equivariant bounded operator
\(\BRA{\xi}\colon \Hils\to\Hilsr\) associated to a
square-integrable vector \(\xi\in\Hils\), where~\(\Hils\) is a
Hilbert space with a coaction of \((\QG,\com)\) and~\(\Hilsr\)
is equipped with the left regular corepresentation of
\((\QG,\com)\).  This construction implies immediately that
square-integrable irreducible corepresentations embed
isometrically and equivariantly into the left regular
corepresentation.

For compact quantum groups, square-integrability is no
restriction, and our theory specialises to the familiar
orthogonality relations for irreducible corepresentations of
compact quantum groups that appear in the Peter--Weyl Theorem.

We study integrable coactions on \(C^*\)\nb-algebras in
Section~\ref{sec:integrable}.  The main ingredient in the
definition is the densely defined, unbounded linear map
\[
\ID_\coef\otimes\Haar\colon
\Mult(\coef\otimes\QG)\supseteq \dom(\ID_\coef\otimes\Haar)
\to \Mult(\coef)
\]
described in~\cite{Kustermans-Vaes:LCQG}.  We call a positive
element \(x\in\coef^+\) integrable with respect to a
coaction~\(\coact_\coef\) of \((\QG,\com)\) on~\(\coef\) if
\(\coact_\coef(x)\) belongs to the domain of
\(\ID_\coef\otimes\Haar\).  The coaction~\(\coact_\coef\) is
integrable if the set of positive integrable elements is dense
in~\(\coef^+\). It follows from the left invariance of the Haar
weight that the range of \(\ID_\coef\otimes\Haar\) consists of
\(\QG\)\nb-invariant multipliers only.  The most important
example of an integrable coaction is the coaction~\(\com\)
of~\(\QG\) on itself.  Using some general permanence properties
of integrability, this example implies that stable coactions
and dual coactions are integrable.  Here stable coactions are
coactions on \(\Comp\bigl(\coef\otimes\Hilsr\bigr)\), where we
equip the Hilbert \(\coef\)\nb-module \(\coef\otimes\Hilsr\)
with the diagonal coaction, using the left regular coaction
on~\(\Hilsr\).

Both square-integrable Hilbert space corepresentations and
integrable coactions on \(C^*\)\nb-algebras are special cases
of square-integrable coactions on Hilbert modules, which we
introduce in Section~\ref{sec:si_Hilm}.  Let~\(\Hilm\) be a
Hilbert \(\coef\)\nb-module with a coaction \(\coact_\Hilm\) of
\((\QG,\com)\).  We define the coefficient
\(c_{\xi\eta}\in\Mult(\coef\otimes\QG)\) for
\(\xi,\eta\in\Hilm\) by
\[
c_{\xi\eta} \defeq \coact_\Hilm(\xi)^*(\eta\otimes1_\QG).
\]
We call~\(\xi\) square-integrable if this belongs to the domain
of the densely defined unbounded map
\(\ID_\coef\otimes\Lambda\colon \Mult(\coef\otimes\QG)\supseteq
\dom(\ID_\coef\otimes\Lambda)\to\coef\otimes\Hilsr\) for all
\(\eta\in\Hilm\); the latter map is also described
in~\cite{Kustermans-Vaes:LCQG}.  The Hilbert module~\(\Hilm\)
is called square-integrable if square-integrable elements are
dense.

We check that \(\xi\in\Hilm\) is square-integrable if and only
if the compact operator \(\ket{\xi}\bra{\xi}\) is integrable
with respect to the coaction on \(\Comp(\Hilm)\) induced by the
coaction on~\(\Hilm\).  This allows us to carry over many
results on integrable coactions on \(C^*\)\nb-algebras.  In
particular, we conclude that \(\Hilm\otimes\Hilsr\) with the
diagonal coaction is square-integrable for any~\(\Hilm\).  We
also establish that the operator \(\BRA{\xi}\colon\Hilm\to
\coef\otimes\Hilsr\) that is defined for
square-integrable~\(\xi\) by \(\BRA{\xi}(\eta) =
(\ID\otimes\Lambda)(c_{\xi\eta})\) is adjointable and
\(\QG\)\nb-equivariant, and we describe its adjoint.

As a consequence, a square-integrable Hilbert module admits
many equivariant adjointable maps to the standard Hilbert
module \(\coef\otimes\Hilsr\).  This together with an idea by
Mingo and Phillips (see~\cite{Mingo-Phillips:Triviality}) is
the main ingredient in our proof of the equivariant
stabilisation theorem, which occupies
Section~\ref{sec:stabilisation}.  It asserts that a countably
generated \(\QG\)\nb-equivariant Hilbert
\(\coef\)\nb-module~\(\Hilm\) is square-integrable if and only
if there is an equivariant isomorphism
\[
\Hilm\oplus \coef\otimes\Hilsr^\infty \cong
\coef\otimes\Hilsr^\infty.
\]
Here \(\coef\otimes\Hilsr\) carries the diagonal coaction using
the left regular corepresentation on~\(\Hilsr\).  Thus
\(\coef\otimes\Hilsr\) plays the role of \(L^2(G,\coef)\) in
the group case.

\section{Preliminaries on locally compact quantum groups}
\label{sec:LCQG}

We use the theory of locally compact quantum groups developed
by Johan Kustermans and Stefaan Vaes in
\cites{Kustermans-Vaes:LCQG,Kustermans-Vaes:LCQGvN} throughout
this article.  In this section, we fix our notation and recall
some important facts that we shall need repeatedly.

The theory of square-integrable \emph{group} representations is
based on a close relationship between the spaces
\(\CONT_0(G)\), \(\Cred(G)\), and \(L^2(G)\).  Let \(f\colon
G\to\C\) be a continuous function with compact support.  Such a
function can play at least three different roles.  First, we
may view~\(f\) as an element of the \(C^*\)\nb-algebra
\(\CONT_0(G)\) of continuous functions on~\(G\) vanishing at
infinity.  Secondly, we may view~\(f\) as an element of the
reduced group \(C^*\)\nb-algebra \(\Cred(G)\).  And thirdly, we
may view~\(f\) as an element of the Hilbert space \(L^2(G) =
L^2(G,\dd x)\) of square-integrable functions on~\(G\); here
\(\dd x\) denotes a left invariant Haar measure on~\(G\).  More
formally, the comparison of these three different roles
provides us with densely defined, closed, unbounded operators
between the Banach spaces \(\CONT_0(G)\), \(L^2(G)\), and
\(\Cred(G)\).  We shall need a quantum group generalisation of
these unbounded operators.

Let~\(\QG\) be a locally compact quantum group and
let~\(\hat{\QG}\) be its dual.  Their quantum group structures
are encoded by comultiplications
\[
\com\colon \QG\to\Mult(\QG\otimes\QG),
\qquad
\hat{\com}\colon \hat{\QG}\to
\Mult\bigl(\hat{\QG}\otimes\hat{\QG}\bigr),
\]
where~\(\Mult\) denotes multiplier algebras and~\(\otimes\)
spatial \(C^*\)\nb-algebra tensor products.  A basic
requirement in~\cite{Kustermans-Vaes:LCQG} is the existence of
left and right invariant faithful proper weights on~\(\QG\)
which satisfy a KMS condition; the existence of such weights
for the dual follows.  We denote the left invariant Haar
weights on \(\QG\) and~\(\hat{\QG}\) by \(\Haar\)
and~\(\hat{\Haar}\).

Any proper weight on a \(C^*\)\nb-algebra gives rise to a
Hilbert space representation by a generalisation of the
Gelfand--Naimark--Segal construction (see~\cite{Combes:Poids}).
When we apply this to the weight~\(\Haar\), we get a Hilbert
space~\(\Hilsr\), a faithful \Star{}representation of~\(\QG\)
on~\(\Hilsr\), which we simply denote by left multiplication
\(\QG\times\Hilsr\to\Hilsr\), \((x,\eta)\mapsto x\cdot\eta\),
and a closed unbounded linear map~\(\Lambda\) between dense
subspaces of~\(\QG\) and~\(\Hilsr\).

Actually, we mainly use the extension of~\(\Lambda\) to the
multiplier algebra.  The weight~\(\Haar\) extends uniquely to a
strictly lower semi-continuous weight on \(\Mult(\QG)\), which
we still denote by~\(\Haar\).  We let
\begin{equation}
  \label{eq:dom_Lambda}
  \dom \Lambda
  \defeq \{x\in\Mult(\QG)\mid \Haar(x^*x)<\infty\}.
\end{equation}
There is a linear map \(\Lambda\colon \dom\Lambda\to\Hilsr\)
that satisfies
\begin{alignat}{2}
  \label{eq:inner_product_Lambda}
  \bigl\langle \Lambda(x),\Lambda(y) \bigr\rangle &=
  \Haar(x^*y)
  &\qquad&\text{for all \(x,y\in\dom\Lambda\), and}\\
  \label{eq:Lambda_pi}
  \Lambda(xy) &= x\cdot\Lambda(y)
  &\qquad&\text{for all \(x\in\Mult(\QG)\),
    \(y\in\dom \Lambda\).}
\end{alignat}
The map~\(\Lambda\) above is closed, densely defined, and has
dense range with respect to the strict topology on
\(\Mult(\QG)\) and the norm topology on~\(\Hilsr\).  Its
restriction to~\(\QG\) itself is closed, densely defined, and
has dense range with respect to the norm topologies on both
\(\QG\) and~\(\Hilsr\).

Notice that unlike in
\cites{Kustermans-Vaes:LCQG,Kustermans-Vaes:LCQGvN} our Hilbert
space inner products are linear in the second variable---this
is standard for Hilbert modules.

In the group case, we have \(\QG=\CONT_0(G)\) and
\(\Hilsr=L^2(G)\), where we tacitly use the left invariant Haar
measure on~\(G\).  Then \(\dom\Lambda = \CONT_0(G)\cap
L^2(G)\), viewed as a subspace of \(\CONT_0(G)\),
and~\(\Lambda\) maps this “identically” to a subspace of
\(L^2(G)\).

If~\(G\) is compact, then~\(\Lambda\) is everywhere defined and
bounded, that is, \(\CONT_0(G)=\CONT(G)\) is contained
in~\(L^2(G)\).  This observation generalises to a
characterisation for compact quantum groups:

\begin{remark}
  \label{rem:compact_qg}
  A quantum group~\(\QG\) is compact if and only if the Haar
  weight~\(\Haar\) is bounded, if and only if~\(\Lambda\) is a
  bounded linear map \(\QG\to\Hilsr\).
\end{remark}

The quantum group structure on~\(\QG\) may also be encoded
using a \emph{multiplicative unitary}.  This is the unique
unitary operator~\(W\) on the Hilbert space tensor product
\(\Hilsr\otimes\Hilsr\) that satisfies
\begin{equation}
  \label{eq:mult_unitary}
  W^*\bigl(\Lambda(x)\otimes\Lambda(y)\bigr) =
  (\Lambda\otimes\Lambda)
  \bigl(\com(y)\cdot (x\otimes 1)\bigr)
  \qquad\text{for all \(x,y\in\dom\Lambda\).}
\end{equation}
Then
\begin{equation}
  \label{eq:comul_via_mult_unitary}
  \com(x) = W^*(1\otimes x)W
  \qquad\text{for all \(x\in\QG\),}
\end{equation}
where we view both sides as operators on~\(\Hilsr\).

Let~\(\QG^\dual\) be the dual Banach space of bounded linear
functionals on~\(\QG\).  Let \(\Bound(\Hilsr)\) be the algebra
of bounded operators on~\(\Hilsr\).  We define a map
\[
\lambda\colon \QG^\dual\to\Bound(\Hilsr),
\qquad
\omega\mapsto (\omega\otimes\ID)(W)
\]
as in \cite{Kustermans-Vaes:LCQGvN}*{\S1.1}.  This becomes an
algebra homomorphism if we equip~\(\QG^\dual\) with the
convolution product \(\omega\star\eta\defeq
(\omega\otimes\eta)\circ\com\) for
\(\omega,\eta\in\hat{\QG}\); here we have tacitly extended
\(\omega\otimes\eta\) to a strictly continuous linear
functional on \(\Mult(\QG\otimes\QG)\).

In the group case, \(\QG^\dual\) is the Banach space of bounded
measures on~\(G\); the map~\(\lambda\) lets a measure act on
\(L^2(G)\) by convolution on the left; and~\(\star\) is the
usual convolution of measures.  To get \(\Cred(G)\), we must
restrict attention to the subspace \(L^1(G)\) of measures that
are absolutely continuous with respect to the Haar measure.
Something similar happens in the quantum group case.

If \(a,b\in\QG\), then we define possibly unbounded linear
functionals \(b\cdot\Haar\), \(\Haar\cdot a\), and
\(b\cdot\Haar\cdot a\) on~\(\QG\) by
\[
b\cdot\Haar(x) \defeq \Haar(x b),\qquad
\Haar\cdot a(x) \defeq \Haar(a x),\qquad
b\cdot\Haar\cdot a(x) \defeq \Haar(a x b).
\]
The functional \(b\cdot\Haar\cdot a\) is bounded if
\(a^*,b\in\dom\Lambda\) because
\begin{equation}
  \label{eq:Haar_as_inner_product}
  b\cdot\Haar\cdot a(x)
  = \bigl\langle \Lambda(a^*),x\cdot\Lambda(b) \bigr\rangle
  \qquad\text{for all \(a^*,b\in\dom\Lambda\), \(x\in\QG\).}
\end{equation}
Let \(L^1(\QG)\subseteq \QG^\dual\) be the closed
linear span of the set of functionals \(b\cdot\Haar\cdot a^*\)
for \(a,b\in\dom\Lambda\).

We define an unbounded linear operator~\(\QGtodual\)
from~\(\QG\) to \(L^1(\QG)\subseteq\QG^\dual\) by letting
\begin{equation}
  \label{eq:QGtodual_dom}
  \begin{aligned}
    x&\in \dom \QGtodual \iff
    \text{\(x\cdot\Haar\)
      is bounded and belongs to
      \(L^1(\QG)\subseteq\QG^\dual\)}\\
    \QGtodual(x) &\defeq x\cdot\Haar
    \qquad\text{for all \(x\in\dom\QGtodual\).}
  \end{aligned}
\end{equation}
It is easy to see that this operator is closed.  We will see
below that it is densely defined and has dense range.

The \(C^*\)\nb-algebra~\(\hat{\QG}\) is, by definition, the
closure of \(\lambda\bigl(L^1(\QG)\bigr)\) in
\(\Bound(\Hilsr)\).  Its comultiplication is defined by
\begin{equation}
  \label{eq:comul_dual}
  \hat{\com}(x)
  \defeq \hat{W}^* (1\otimes x)\hat{W}
  \qquad\text{with}\quad \hat{W}\defeq \Sigma W^*\Sigma,
\end{equation}
where~\(\Sigma\) flips the two tensor factors in
\(\Hilsr\otimes\Hilsr\).

\begin{remark}
  \label{rem:vector_functional}
  Since the von Neumann algebra generated by~\(\QG\) in
  \(\Bound(\Hilsr)\) is in standard form
  (see~\cite{Stratila-Zsido:von_Neumann}*{10.15}), all normal
  functionals on it are vector functionals by
  \cite{Takesaki:Theory_1}*{Theorem V.3.15}).  Therefore, any
  element of \(L^1(\QG)\) is of the form \(x\mapsto
  \langle \eta, x\omega\rangle\) for some
  \(\eta,\omega\in\Hilsr\).  This generalises the fact that any
  function in \(L^1(G)\) is a product of two functions in
  \(L^2(G)\).
\end{remark}

Recall that~\(\Haar\) is a KMS weight with respect to the
\emph{modular automorphism group} \(\modauto\colon
\R\times\QG\to\QG\), which is determined uniquely by the
weight~\(\Haar\).  Extending~\(\modauto\) analytically
from~\(\R\) to~\(\C\), the KMS condition implies
\begin{equation}
  \label{eq:KMS_1}
  \Haar(yx) = \Haar\bigl(x\modauto_{-\ima}(y)\bigr)
  \qquad\text{for all \(y\in\dom\modauto_{-\ima}\),
    \(x\in\dom\Haar\)}
\end{equation}
(see~\cite{Kustermans-Vaes:LCQG}*{Proposition 1.12.(3)}).
Briefly, we have
\begin{equation}
  \label{eq:KMS_2}
  y\cdot\Haar = \Haar\cdot\modauto_{-\ima}(y)
  \qquad\text{for all \(y\in\dom\modauto_{-\ima}\).}
\end{equation}
It follows that the functionals \(b\cdot\Haar\) and
\(\Haar\cdot a\), and \(b\cdot\Haar\cdot a\) are densely
defined for all \(a,b\in\QG\).  Moreover, \(ab\cdot\Haar(x) =
\bigl(b\cdot\Haar\cdot\modauto_{-\ima}(a)\bigr)(x) =
\bigl\langle \Lambda(b^*), x\cdot
\Lambda(\modauto_{-\ima}a)\bigr\rangle\), so that
\(ab\cdot\Haar\) belongs to
\(L^1(\QG)\subseteq\hat{\QG}\) if \(b^*\in\dom\Lambda\)
and \(a\in\dom (\Lambda\circ\modauto_{-\ima})\).

The subsets of \(ab\cdot\Haar\) and of \(ab\) with
\(a,b\) as above are dense in \(L^1(\QG)\) and
in~\(\QG\), respectively.  This shows that the operator
\(\QGtodual\colon \QG\supseteq \dom\QGtodual\to
L^1(\QG)\) above is densely defined and has dense range.

We also get a densely defined unbounded linear
map~\(\hat{\Lambda}\) from~\(\hat{\QG}\) to~\(\Hilsr\) by
\[
\hat{\Lambda}\bigl(\lambda(x\cdot\Haar)\bigr) \defeq \Lambda(x)
\qquad
\text{for \(x\in \dom\QGtodual\)}
\]
(we have \(\dom\QGtodual\subseteq \dom\Lambda\)).  More
precisely, we let~\(\hat{\Lambda}\) be the closure of this
linear map, which is a densely defined closed linear operator
from~\(\hat{\QG}\) to~\(\Hilsr\).

This definition of~\(\hat{\Lambda}\) is equivalent to the one
in \cite{Kustermans-Vaes:LCQGvN}*{page 75}.  As a consequence,
\(\hat{\Lambda}\) together with the identical representation
of~\(\hat{\QG}\) on~\(\Hilsr\) is a GNS construction for the
Haar weight~\(\hat{\Haar}\) on~\(\hat{\QG}\), that is,
\begin{equation}
  \label{eq:dual_Haar}
  \hat{\Haar}(x^*y)
  = \bigl\langle\hat{\Lambda}(x),
  \hat{\Lambda}(y)\bigr\rangle,\qquad
  \hat{\Lambda}(ay) = a\bigl(\hat{\Lambda}(y)\bigr)
\end{equation}
for all \(x,y\in\dom\hat{\Lambda}\),
\(a\in\hat{\QG}\subseteq\Bound(\Hilsr)\).  In particular,
\begin{equation}
  \label{eq:dual_Haar_2}
  \hat{\Haar}\bigl( \lambda(x\cdot\Haar)^*\cdot
  \lambda(y\cdot\Haar)\bigr)
  = \bigl\langle \Lambda(x),\Lambda(y)\bigr\rangle
  = \Haar(x^*y)
\end{equation}
for all \(x,y\in\dom\QGtodual\).

In the group case, \(x\mapsto x\cdot\Haar\) is the standard way
to map a function on~\(G\) to a measure on~\(G\).  Both
\(\QGtodual\) and~\(\hat{\Lambda}\) extend “identical” maps on
suitable dense subspaces: \(\QGtodual\) is the identical map
from \(\CONT_0(G)\cap L^1(G)\) viewed as a subspace of
\(\CONT_0(G)\) to the same space viewed as a subspace of
\(L^1(G)\), and~\(\hat{\Lambda}\) extends the identical map
from \(L^1(G)\cap L^2(G)\) viewed as a subspace of \(\Cred(G)\)
to the same space viewed as a subspace of \(L^2(G)\).

\begin{remark}
  \label{rem:discrete_qg_lambda}
  Let~\(\QG\) be compact, so that~\(\Lambda\) is bounded.  If
  \(x\in \dom\QGtodual\), then
  \begin{multline*}
    \norm{x\Haar}
    = \sup \bigl\{\abs{x\Haar(y^*)} \bigm| y\in\QG,\
    \norm{y}\le1 \bigr\}
    = \sup \bigl\{\abs{\langle \Lambda(y),\Lambda(x)\rangle}
    \bigm| y\in\QG,\ \norm{y}\le1\bigr\}
    \\\le \norm{\Lambda}\cdot \sup \bigl\{\abs{\langle
      \Lambda(y),\Lambda(x)\rangle} \bigm|
    y\in\QG,\ \norm{\Lambda(y)}\le1\bigr\}
    =  \norm{\Lambda}\cdot \norm{\Lambda(x)},
  \end{multline*}
  where we used the definition of the norm on~\(\QG^\dual\),
  the relation \(x\Haar(y^*)= \bigl\langle
  \Lambda(y),\Lambda(x) \bigr\rangle\), and Riesz' Theorem for
  the Hilbert space~\(\Hilsr\).  Since
  \(\Lambda(\dom\QGtodual)\) is dense in~\(\Hilsr\), it follows
  that there is a bounded map \(\Hilsr\to L^1(\QG)\) that maps
  \(\hat{\Lambda}\bigl(\lambda(\omega)\bigr)\) to~\(\omega\)
  for all \(\omega\in L^1(\QG)\).

  Since any discrete quantum group is the dual of a compact
  one, this shows that the map \(\Lambda\colon \QG\supseteq
  \dom\Lambda\to\Hilsr\) is the inverse of a bounded map
  \(\Lambda^{-1}\colon \Hilsr\to \QG\) for any discrete quantum
  group~\(\QG\).
\end{remark}

Recall that a (right) \emph{corepresentation} of~\(\QG\) on a
Hilbert space~\(\Hils\) is a unitary operator~\(U\) on the
Hilbert \(\QG\)\nb-module \(\Hils\otimes\QG\) that satisfies
\((\ID_\Hils\otimes\com)(U) = U_{12}\cdot U_{13}\), where we
use the standard leg numbering notation.  This generates a
(right) coaction
\begin{equation}\label{eq:coact_Hils}
\coact_U\colon \Hils
\to \Mult(\Hils\otimes\QG)\defeq \Bound(\QG,\Hils\otimes\QG),
\qquad
\eta\mapsto U(\eta\otimes1_\QG).
\end{equation}
Conversely, any coaction on~\(\Hils\) with the usual properties
comes from a unique corepresentation~\(U\) (see
\cite{Baaj-Skandalis:Hopf_KK}*{Proposition 2.4}).  We will
mainly use right corepresentations and right coactions in this
article.

A right corepresentation of~\(\QG\) yields a \emph{left}
\(\QG^\dual\)\nb-module structure via
\begin{equation}
  \label{eq:integrate_rep}
  \omega \star_U \eta \defeq
  (\ID\otimes\omega)\bigl(\coact_U(\eta)\bigr)
  = (\ID\otimes\omega)\bigl(U(\eta\otimes1)\bigr)
  \eqdef (\ID\otimes\omega)(U)(\eta).
\end{equation}
This module structure is always non-degenerate.  The
\emph{universal locally compact quantum group}
\((\hat{\QG}_\uni,\hat{\com}_\uni)\) associated to the reduced
locally compact quantum group \((\hat{\QG},\hat{\com})\) is
defined by the universal property that its non-degenerate
\Star{}representations correspond exactly to right
corepresentations of~\(\QG\),
see~\cite{Kustermans:LCQG_universal}.

In the group case, a strongly continuous unitary group
representation~\(\pi\) of~\(G\) on~\(\Hils\) is encoded either
by the corepresentation
\[
U\colon \CONT_0(G,\Hils)\to\CONT_0(G,\Hils),
\qquad (Uf)(g) \defeq \pi_g\bigl(f(g)\bigr),
\]
or by the coaction
\[
\coact_U\colon \Hils\to \CONT_b(G,\Hils),
\qquad (\coact_U\eta)(g) \defeq \pi_g(\eta).
\]
The integrated form defined above lets a measure~\(\omega\) act
by \(\omega\star\eta = \int_G \pi_g(\eta)\,\dd\omega(g)\).  The
universal dual quantum group~\(\hat{\QG}_\uni\) is the full
group \(C^*\)\nb-algebra, whereas~\(\hat{\QG}\) itself is the
reduced group \(C^*\)\nb-algebra of~\(G\).

The \emph{left regular corepresentation} of~\(\QG\)
on~\(\Hilsr\) is a crucial example for us: it is the paradigm
of square-integrability.  As a right corepresentation, it is
the unitary \(\Sigma W\Sigma\), where~\(W\) is the
multiplicative unitary of~\(\QG\) described
in~\eqref{eq:mult_unitary}, but viewed as an adjointable
operator on the Hilbert \(\QG\)\nb-module \(\QG\otimes\Hilsr\),
and~\(\Sigma\) is the isomorphism
\(\Hilsr\otimes\QG\leftrightarrow\QG\otimes\Hilsr\) that flips
the tensor factors.  The corresponding right
coaction~\(\coact_\lambda\) is
\begin{equation}
  \label{eq:coact_lambda}
  \coact_\lambda(\eta) = \Sigma W(1\otimes\eta)
  \qquad\text{for all \(\eta\in\Hilsr\).}
\end{equation}

Let \(\omega\star_\lambda x\) for \(\omega\in\QG^\dual\),
\(x\in\Hilsr\) denote the resulting non-degenerate left
\(\QG^\dual\)\nb-module structure on~\(\Hilsr\).  The module
structure~\(\star_\lambda\) determines a \Star{}representation
of~\(\hat{\QG}\), not just of~\(\hat{\QG}_\uni\).  This is the
left regular representation~\(\lambda\) we used to
define~\(\hat{\QG}\), that is,
\begin{equation}
  \label{eq:integrated_left_regular}
  \omega \star_\lambda \eta = \lambda(\omega)(\eta)
  \qquad\text{for all \(\omega\in\QG^\dual\),
    \(\eta\in\Hilsr\).}
\end{equation}
This follows from the computation
\begin{multline*}
  \omega \star_\lambda \eta =
  (\ID\otimes\omega)\circ\Sigma\circ W(1\otimes\eta)
  = (\omega\otimes\ID)\bigl(W(1\otimes\eta)\bigr)
  \\= \bigl((\omega\otimes\ID)(W)\bigr)(\eta)
  = \lambda(\omega)(\eta).
\end{multline*}

We return to the regular representation.  The \emph{left
  regular representation}~\(\lambda\) of the dual quantum
group~\(\hat{\QG}\) on~\(\Hilsr\) is described either by
\(\lambda(\omega) = (\omega\otimes\ID)(W)\) as above or by
\(\lambda(\omega)\bigl(\hat{\Lambda}(x)\bigr) =
\hat{\Lambda}(\omega x)\) for all \(\omega\in\hat{\QG}\),
\(x\in \dom\hat{\Lambda}\).  Let~\(\hat{T}\) be the closure of
the unbounded conjugate-linear operator on~\(\Hilsr\) defined
by \(\hat{T}\bigl(\hat{\Lambda}(x)\bigr) \defeq
\hat{\Lambda}(x^*)\) if \(x,x^*\in\dom\hat{\Lambda}\).  This
operator admits a polar decomposition
\begin{equation}
  \label{eq:polar_decompose_adjoint}
  \hat{T} = \hat{J}\hat{\nabla}^{\nicefrac{1}{2}}
  = \hat{\nabla}^{-\nicefrac{1}{2}}\hat{J}
\end{equation}
for some conjugate-unitary~\(\hat{J}\), called \emph{modular
  conjugation}, and a positive unbounded
operator~\(\hat{\nabla}\) on~\(\Hilsr\), called the
\emph{modular operator} of the KMS weight~\(\hat{\Haar}\).
This operator is closely related to the \emph{modular
  automorphism group}~\(\hat{\modauto}\) of~\(\hat{\QG}\):
\begin{equation}
  \label{eq:nabla_generates_modauto}
  \hat{\nabla}^{\ima s} \lambda(\omega) \hat{\nabla}^{-\ima s}
  = \lambda\bigl(\hat{\modauto}_s(\omega)\bigr)
\end{equation}
for all \(s\in\C\), \(\omega\in\dom \hat{\modauto}_s\).

The \emph{right regular anti-representation} of~\(\hat{\QG}\)
is the map
\begin{equation}
  \label{eq:right_regular_anti}
  \rho\colon \hat{\QG} \to \Bound(\Hilsr),
  \qquad \omega\mapsto \hat{J} \lambda(\omega)^*\hat{J},
\end{equation}
which satisfies \(\rho(x^*)=\rho(x)^*\) and
\(\rho(xy)=\rho(y)\rho(x)\) for all \(x,y\in\hat{\QG}\).  It
can be turned into a representation using the unitary antipode
of~\(\hat{\QG}\).  Its range \(\hat{\QG}^\comm \defeq
\hat{J}\hat{\QG}\hat{J}\) is the \emph{\(C^*\)\nb-commutant}
of~\(\hat{\QG}\), which is a locally compact quantum group in
its own right (see~\cite{Kustermans-Vaes:LCQGvN}).  We use
\eqref{eq:nabla_generates_modauto}
and~\eqref{eq:polar_decompose_adjoint} to compute
\begin{multline*}
  \rho(\omega^*) \hat{\Lambda}(x)
  = \hat{J} \lambda(\omega) \hat{J}\hat{\Lambda}(x)
  = \hat{J} \hat{\nabla}^{\nicefrac{1}{2}}
  \lambda\bigl(\hat{\modauto}_{\ima/2}(\omega)\bigr)
  \hat{\nabla}^{-\nicefrac{1}{2}} \hat{J}\hat{\Lambda}(x)
  \\= \hat{T} \lambda\bigl(\hat{\modauto}_{\ima/2}(\omega)\bigr)
  \hat{T}\hat{\Lambda}(x)
  = \hat{\Lambda}(x\cdot \hat{\modauto}_{\ima/2}(\omega)^*)
  = \hat{\Lambda}\bigl(x\cdot
  \hat{\modauto}_{-\ima/2}(\omega^*)\bigr)
\end{multline*}
provided \(x,x^*\in\dom\hat{\Lambda}\) and \(\omega\in\dom
\hat{\modauto}_{\ima/2}\) (then \(x\cdot
\hat{\modauto}_{-\ima/2}(\omega^*) \in \dom\hat{\Lambda}\)).
As a result,
\begin{equation}
  \label{eq:right_regular_rep}
  \rho\bigl(\hat{\modauto}_{\ima/2}(\omega)\bigr)
  \hat{\Lambda}(x)
  = \hat{\Lambda}(x\cdot\omega)
  \qquad\text{for all \(\omega\in \dom
    \hat{\modauto}_{\ima/2}\), \(x\in\dom\hat{\Lambda}\).}
\end{equation}
See also \cite{Kustermans-Vaes:LCQG}*{Proposition 1.12.(2)}.

\section{Square-integrable Hilbert space representations}
\label{sec:si_Hilbert_space}

In the Hilbert space case, it makes essentially no difference
whether we study corepresentations of~\(\QG\) or
representations of \(\hat{\QG}_\uni\).  The latter point of
view has the advantage that square-integrable representations
in this context are a special case of a general notion due to
François Combes~\cite{Combes:Systemes_Hilbertiens}.  Here we
recall this definition, check that the left regular
representation is square-integrable, and establish the
orthogonality relations for coefficients of square-integrable
representations.  Our proofs follow those by Marc Rieffel
in~\cite{Rieffel:Integrable_proper}.

The universal dual~\(\hat{\QG}_\uni\) has its own Haar weight
and corresponding GNS construction.  Actually, we get these by
composing the corresponding maps for~\(\hat{\QG}\) with the
canonical quotient mapping \(\hat{\QG}_\uni\prto\hat{\QG}\)
(see~\cite{Kustermans:LCQG_universal}).  Therefore, we denote
the Haar weight on~\(\hat{\QG}_\uni\) by~\(\hat{\Haar}\) as
well, and we still write~\(\hat{\Lambda}\) for the unbounded
map from~\(\hat{\QG}_\uni\) to the Hilbert space~\(\Hilsr\).
The \(C^*\)\nb-algebra~\(\hat{\QG}_\uni\) with the left ideal
\(\dom \hat{\Lambda}\) and the sesquilinear form \((x,y)\mapsto
\hat{\Haar}(y^*x)\) is a left Hilbert system in the notation
of~\cite{Combes:Systemes_Hilbertiens}.

Let \(U\in\Mult(\Hils\otimes\QG)\) be a corepresentation
of~\(\QG\) on a Hilbert space~\(\Hils\), and let \(\rho_U\colon
\hat{\QG}_\uni \to \Bound(\Hils)\) be the corresponding
\Star{}representation of~\(\hat{\QG}_\uni\).

\begin{definition}[\cite{Combes:Systemes_Hilbertiens}*{\S1.4}]
  \label{def:coefficient_Hilbert}
  Let \(\alpha,\beta\in\Hils\).  The \emph{coefficient}
  \(c_{\alpha\beta}\) of the representation~\(\rho_U\) is the
  linear functional
  \[
  \hat{\QG}_\uni\to\C,
  \qquad x\mapsto \langle \beta,\rho_U(x)\alpha\rangle.
  \]
  (Notice that our inner products are linear in the second
  variable.)

  This coefficient is called \emph{square-integrable} if there
  is \(r\in\R_{>0}\) with
  \[
  \abs{c_{\alpha\beta}(x)}^2
  \le r\cdot \hat{\Haar}(x^*x)
  = r\cdot \norm{\hat{\Lambda}(x)}_{\Hilsr}^2.
  \]
  In this case, \(\hat{\Lambda}(x)\mapsto
  c_{\alpha\beta}(x)\) extends uniquely to a bounded linear
  functional on~\(\Hilsr\).  This is of the form
  \(c_{\alpha\beta}(x) = \bigl\langle \tilde{c}_{\alpha\beta},
  \hat{\Lambda}(x)\bigr\rangle\) for some
  \(\tilde{c}_{\alpha\beta}\in\Hilsr\) by Riesz's Theorem.
\end{definition}

\begin{definition}
  \label{def:si_Hilbert}
  We call~\(\alpha\) \emph{square-integrable} (with respect
  to~\(U\)) if~\(c_{\alpha\beta}\) is square-integrable for
  each \(\beta\in\Hils\).  We call~\(\alpha\)
  \emph{\(U\)\nb-bounded} if if there is \(r\in\R_{>0}\) with
  \[
  \norm{\rho_U(x)(\alpha)}^2
  \le r\cdot \norm{\hat{\Lambda}(x)}_{\Hilsr}^2
  = r\cdot \hat{\Haar}(x^*x)
  \qquad\text{for all \(x\in\dom\hat{\Lambda}\).}
  \]
  We let \(\Hils_\si\subseteq\Hils\) be the subset of
  square-integrable vectors.

  The corepresentation~\(U\) is called \emph{square-integrable}
  if the subspace~\(\Hils_\si\) of square-integrable vectors is
  dense in~\(\Hils\).
\end{definition}

\begin{lemma}[\cite{Combes:Systemes_Hilbertiens}*{\S1.4}]
  \label{lem:bounded_si}
  A vector is \(U\)\nb-bounded if and only
  if it is square-integrable.
\end{lemma}

\begin{proof}
  It is clear that bounded vectors are square-integrable.

  Conversely, if \(\alpha\in\Hils_\si\), then there is a linear
  map
  \begin{equation}
    \label{eq:def_BRA_Hilbert}
    \BRA{\alpha}\colon \Hils\to\Hilsr,
    \qquad \beta\mapsto \tilde{c}_{\alpha\beta}.
  \end{equation}
  Since its graph is closed, it is bounded.  Let
  \(\KET{\alpha}\defeq \BRA{\alpha}^*\colon \Hilsr\to\Hils\) be
  its adjoint.  We have
  \begin{equation}
    \label{eq:def_KET_Hilbert}
    \KET{\alpha}\bigl(\hat{\Lambda}(x)\bigr) =
    \rho_U(x)(\alpha)
    \qquad\text{for all \(x\in\dom\hat{\Lambda}\)}
  \end{equation}
  because \(\bigl\langle \beta, \rho_U(x)(\alpha)\bigr\rangle =
  c_{\alpha\beta}(x) = \bigl\langle \tilde{c}_{\alpha\beta},
  \hat{\Lambda}(x)\bigr\rangle = \bigl\langle
  \BRA{\alpha}\beta, \hat{\Lambda}(x)\bigr\rangle =
  \bigl\langle \beta,
  \KET{\alpha}\hat{\Lambda}(x)\bigr\rangle\).  The boundedness
  of~\(\KET{\alpha}\) means that~\(\alpha\) is
  \(U\)\nb-bounded.
\end{proof}

The operators \(\KET{\alpha}\) and \(\BRA{\alpha}\) introduced
in the proof of Lemma~\ref{lem:bounded_si} are the cornerstones
of the general theory of square-integrable representations.

Since \(\hat{\Lambda}\) factors through the projection
\(\hat{\QG}_\uni\to\hat{\QG}\), Lemma~\ref{lem:bounded_si}
implies that any square-integrable representation
of~\(\hat{\QG}_\uni\) factors through
\(\hat{\QG}_\uni\to\hat{\QG}\), so that it is already a
representation of~\(\hat{\QG}\).  Therefore, we may
replace~\(\hat{\QG}_\uni\) by~\(\hat{\QG}\) in the following.

\begin{lemma}
  \label{lem:KET_equivariant}
  Let \(\alpha\in\Hils_\si\).  The operators
  \(\KET{\alpha}\colon \Hilsr\to\Hils\) and
  \(\BRA{\alpha}\colon \Hils\to\Hilsr\) intertwine the
  corepresentation~\(U\) and the left regular corepresentation
  on~\(\Hilsr\).
\end{lemma}

Recall that an operator \(T\colon \Hils\to\Hilsr\) is called an
\emph{intertwining operator} for the coactions
\(\coact_U\) and~\(\coact_\lambda\) if \(\coact_\lambda(Tx) =
(T\otimes\ID)\circ\coact_U(x)\) for all \(x\in\Hils\).

\begin{proof}
  The operator \(\KET{\alpha}\) is a left
  \(\hat{\QG}\)\nb-module homomorphism
  by~\eqref{eq:def_KET_Hilbert}.  So is its
  adjoint~\(\BRA{\alpha}\) because the representations
  of~\(\hat{\QG}\) on \(\Hils\) and~\(\Hilsr\) are involutive.
  But the \(\hat{\QG}\)\nb-module homomorphisms with respect to
  the integrated forms of two corepresentations are exactly the
  intertwining operators.
\end{proof}

\begin{lemma}
  \label{lem:regular_si}
  The left regular corepresentation of~\(\QG\) on~\(\Hilsr\) is
  square-integrable.
\end{lemma}

\begin{proof}
  Let~\(\hat{\modauto}\) denote the modular automorphism group
  of the weight~\(\hat{\Haar}\) on~\(\hat{\QG}\) and
  let~\(\rho\) be the right regular anti-representation
  of~\(\hat{\QG}\) defined in~\eqref{eq:right_regular_anti}.
  Equation~\eqref{eq:right_regular_rep} yields
  \[
  \KET{\hat{\Lambda}(x)}\hat{\Lambda}(a)
  = \hat{\Lambda}(ax)
  = \rho\bigl(\hat{\modauto}_{\ima/2}(x)\bigr)
  \hat{\Lambda}(a)
  \]
  for all \(a\in\dom\hat{\Lambda}\), \(x\in
  \dom\hat{\Lambda}\cap \dom \hat{\modauto}_{\ima/2}\).
  Since~\(\rho\) is an isometric \Star{}anti-representation
  of~\(\hat{\QG}\), this shows that \(\bigl\lVert
  \KET{\hat{\Lambda}(x)}\bigr\rVert =
  \norm{\hat{\modauto}_{\ima/2}(x)}<\infty\) for all \(x\in
  \dom\hat{\Lambda}\cap \dom \hat{\modauto}_{\ima/2}\), so that
  these vectors are \(\lambda\)\nb-bounded and hence
  square-integrable.  Since~\(\hat{\Lambda}\) maps
  \(\dom\hat{\Lambda}\cap \dom \hat{\modauto}_{\ima/2}\) to a
  dense subspace of~\(\Hilsr\), the square-integrable vectors
  are dense for the left regular corepresentation.
\end{proof}

\begin{theorem}
  \label{the:si_stabilisation_Hils}
  A Hilbert space corepresentation is square-integrable if and
  only if it is a direct sum of sub-corepresentations of the
  left regular corepresentation, if and only if its integrated
  form \(\hat{\QG}_\uni\to\Bound(\Hils)\) extends to a normal
  \Star{}representation of the von Neumann algebra completion
  \(\lambda(\hat{\QG})''\) of~\(\hat{\QG}\).
\end{theorem}

\begin{proof}
  It is easy to see that square-integrability is inherited by
  sub-corepresentations and direct sums.  Hence
  Lemma~\ref{lem:regular_si} shows that direct sums of
  sub-corepresentations of the left regular corepresentation
  are square-integrable.

  Conversely, let~\(U\) be a square-integrable corepresentation
  on a Hilbert space~\(\Hils\).  Any \(\xi\in\Hils_\si\) yields
  a non-zero equivariant linear operator \(\KET{\xi}\colon
  \Hilsr\to\Hils\).  The partial isometry from its polar
  decomposition is equivariant as well and identifies a
  non-zero sub-corepresentation of~\(\Hils\) with a
  sub-corepresentation of~\(\Hilsr\).  Now we use Zorn's Lemma
  to complete the proof: there is a maximal set of orthogonal,
  non-zero \(\QG\)\nb-invariant subspaces \((\Hils_i)_{i\in
    I}\) such that the restriction of~\(U\) to each~\(\Hils_i\)
  is unitarily equivalent to a sub-corepresentation of the left
  regular corepresentation.  If \(\bigoplus \Hils_i\) is not
  yet all of~\(\Hils\), then its complement contains a
  square-integrable vector, which yields a sub-corepresentation
  that we may add to our set of subspaces, contradicting
  maximality.  Hence \(\bigoplus \Hils_i=\Hils\).

  Of course, the left regular representation
  of~\(\hat{\QG}_\uni\) on~\(\Hilsr\) extends to a normal
  \Star{}representation of the von Neumann algebra completion
  \(\lambda(\hat{\QG})''\).  This remains true for any direct
  sum of sub-corepresentations of the left regular
  corepresentation and hence holds for any square-integrable
  corepresentation.  Conversely, it is well-known that any such
  normal \Star{}representation of \(\lambda(\hat{\QG})''\) may
  be decomposed as a direct sum of subrepresentations of the
  standard representation.
\end{proof}

\begin{lemma}
  \label{lem:act_on_coefficient_right}
  Let \(\xi\in\Hils_\si\) and let
  \(\omega\in\QG^\dual\subseteq\Mult(\hat{\QG})\).  If
  \(\omega\in\dom\hat{\modauto}_{\nicefrac{\ima}{2}}\), then
  \[
  \rho_U(\omega)(\xi) = \omega\star\xi\in\Hils_\si
  \qquad\text{with \(\KET{\omega\star\xi} = \KET{\xi} \circ
    \rho\bigl(\hat{\modauto}_{\nicefrac{\ima}{2}}(\omega)\bigr)\).}
  \]
\end{lemma}

Here~\(\rho\) denotes the right regular anti-representation
on~\(\Hilsr\) defined in~\eqref{eq:right_regular_anti}.
Lemma~\ref{lem:act_on_coefficient_right} may be proved using
the same ingredients as for Lemma~\ref{lem:regular_si}.  We
omit the details because this is a special case of
Lemma~\ref{lem:si_act} below, anyway.

Now let \(U\) and~\(U'\) be two \emph{irreducible}
corepresentations of~\(\QG\) on Hilbert spaces \(\Hils\)
and~\(\Hils'\).  The \emph{orthogonality relations} for
coefficients of square-integrable representations describe the
inner products \(\langle \tilde{c}_{\alpha\beta},
\tilde{c}_{\gamma\delta}\rangle\) for \(\alpha\in\Hils_\si\),
\(\beta\in\Hils\), \(\gamma\in\Hils'_\si\),
\(\delta\in\Hils'\).  Using the operators just introduced, we
can rewrite this as
\[
\langle \tilde{c}_{\alpha\beta},
\tilde{c}_{\gamma\delta}\rangle
= \bigl\langle \BRA{\alpha}(\beta),
\BRA{\gamma}(\delta)\bigr\rangle
= \bigl\langle \beta,
\KET{\alpha}\BRA{\gamma}(\delta)\bigr\rangle.
\]

\begin{proposition}[First Orthogonality Relation]
  \label{pro:first_orthogonality}
  If \(U\) and~\(U'\) are non-equivalent irreducible
  corepresentations of~\(\QG\), then \(\langle
  \tilde{c}_{\alpha\beta}, \tilde{c}_{\gamma\delta}\rangle=0\)
  for all \(\alpha\in\Hils_\si\), \(\beta\in\Hils\),
  \(\gamma\in\Hils'_\si\), \(\delta\in\Hils'\).
\end{proposition}

\begin{proof}
  The \(\QG\)\nb-equivariant operator
  \(\KET{\alpha}\BRA{\gamma}\colon \Hils'\to\Hils\) vanishes
  because \(U\) and~\(U'\) are not equivalent.
\end{proof}

The second orthogonality relation deals with the more
interesting case \(\Hils=\Hils'\) and \(U=U'\).  In that case,
Schur's Lemma implies that \(\KET{\alpha}\BRA{\gamma}\) is a
scalar multiple of the identity map.  Thus there is a
sesquilinear form~\(s\) on~\(\Hils_\si\) with
\(\KET{\alpha}\BRA{\gamma} = s(\alpha,\gamma)\) for all
\(\alpha,\gamma\in\Hils_\si\).  Since the representation
of~\(\hat{\QG}\) on~\(\Hils\) is non-degenerate,
\(\KET{\alpha}=0\) implies \(\alpha=0\)
by~\eqref{eq:def_KET_Hilbert}.  Hence \(s(\alpha,\alpha)>0\)
for all \(\alpha\in\Hils_\si\) with \(\alpha\neq0\).

Consequently, there is a positive, self-adjoint, unbounded
operator~\(K\) on~\(\Hils\) with \(\ker K=\{0\}\) and
\(s(\alpha,\gamma) = \langle \gamma, K^{-1}\alpha\rangle\) for
all \(\alpha,\gamma\in\Hils_\si\).  That is,
\[
\langle \tilde{c}_{\alpha\beta},
\tilde{c}_{\gamma\delta}\rangle
= \bigl\langle \beta, \langle
\gamma,K^{-1}\alpha\rangle\cdot\delta \bigr\rangle
= \langle \beta,\delta\rangle \cdot \langle
\gamma,K^{-1}\alpha\rangle
= \bigl\langle \beta,\delta\rangle \cdot \langle
K^{-\nicefrac{1}{2}}\gamma,K^{-\nicefrac{1}{2}}\alpha \bigr\rangle
\]
for all \(\alpha,\gamma\in\Hils_\si\),
\(\beta,\delta\in\Hils\).  Taking \(\alpha=\gamma\) and
\(\beta=\delta\), we see that \(\Hils_\si=\dom
K^{-\nicefrac{1}{2}}\).

The positive operator~\(K\) plays the role of the \emph{formal
  dimension} for the square-integrable representation~\(U\),
compare \cite{Rieffel:Integrable_proper}*{\S8} for the case of
group representations.  It remains to describe~\(K\) or,
equivalently, \(K^{-1}\) more explicitly.  This involves the
modular theory of~\(\hat{\QG}\).

Since \(\KET{\xi}\neq0\) for square-integrable \(\xi\neq0\) and
\(\KET{\xi}\BRA{\xi}\) is a scalar multiple of the identity
map, we may choose an auxiliary square-integrable vector
\(\xi\in\Hils_\si\) with \(\KET{\xi}\BRA{\xi}=\ID_\Hils\).
Equivalently, \(\BRA{\xi}\colon \Hils\to\Hilsr\) is an
isometry.

\begin{theorem}[Second Orthogonality Relation]
  \label{the:second_orthogonality}
  Let~\(\QG\) be a locally compact quantum group, let~\(\Hils\)
  be a Hilbert space, and let~\(U\) be an irreducible
  square-integrable corepresentation of~\(\QG\) on~\(\Hils\).
  Let \(\alpha,\gamma\in\Hils_\si\), \(\beta,\delta\in\Hils\),
  and choose \(\xi\in\Hils_\si\) with
  \(\KET{\xi}\BRA{\xi}=\ID_\Hils\).  Let
  \(\hat{\nabla}\colon\Hilsr\supseteq\dom(\hat{\nabla})\to\Hilsr\)
  be the modular operator of the KMS weight~\(\hat{\Haar}\)
  on~\(\hat{\QG}\).  Then
  \[
  \langle \tilde{c}_{\alpha\beta},
  \tilde{c}_{\gamma\delta}\rangle =
  \langle \beta,\delta\rangle \cdot
  \langle K^{-\nicefrac{1}{2}}\gamma,
  K^{-\nicefrac{1}{2}}\alpha\rangle
  \]
  with the positive unbounded operator
  \[
  K^{-1} = \norm{\xi}^{-2}
  \KET{\xi}\circ\hat{\nabla}^{-1}\circ\BRA{\xi}\colon
  \Hils\to\Hils.
  \]
\end{theorem}

The operator~\(K\) is independent of the auxiliary
vector~\(\xi\).

\begin{proof}
  Our main task is to relate the matrix coefficients
  \(\tilde{c}_{\alpha\beta}\) and \(\tilde{c}_{\beta\alpha}\).
  It follows easily from the definition that
  \[
  c_{\alpha\beta}(x)
  = \langle \beta, \rho_U(x)\alpha\rangle
  = \langle \rho_U(x^*)\beta, \alpha\rangle
  = \conj{\langle \alpha, \rho_U(x^*)\beta\rangle}
  = \conj{c_{\beta\alpha}(x^*)}.
  \]
  Equivalently, \(\bigl\langle
  \tilde{c}_{\alpha\beta},\hat{\Lambda}(x) \bigr\rangle =
  \bigl\langle
  \hat{\Lambda}(x^*),\tilde{c}_{\beta\alpha}\bigr\rangle\).

  Let~\(\hat{T}\) be the unbounded, conjugate-linear operator
  on~\(\Hilsr\) that maps \(\hat{\Lambda}(x)\) to
  \(\hat{\Lambda}(x^*)\); recall that it has a polar
  decomposition \(\hat{T} =
  \hat{J}\hat{\nabla}^{\nicefrac{1}{2}} =
  \hat{\nabla}^{-\nicefrac{1}{2}}\hat{J}\) as
  in~\eqref{eq:polar_decompose_adjoint}.

  The computation above shows that
  \begin{equation}
    \label{eq:flip_coefficient}
    \tilde{c}_{\alpha\beta} =
    \hat{T}^*(\tilde{c}_{\beta\alpha})
    \qquad\text{for all \(\alpha,\beta\in\Hils_\si\).}
  \end{equation}
  Hence
  \begin{multline*}
    \langle \xi,\xi\rangle \cdot \langle K^{-\nicefrac{1}{2}}\beta,
    K^{-\nicefrac{1}{2}}\alpha\rangle
    = \langle \tilde{c}_{\alpha\xi}, \tilde{c}_{\beta\xi}\rangle
    = \langle \hat{T}^*\tilde{c}_{\xi\alpha},
    \hat{T}^*\tilde{c}_{\xi\beta}\rangle
    \\= \bigl\langle \hat{J}^* \hat{\nabla}^{-\nicefrac{1}{2}}\BRA{\xi}\alpha,
    \hat{J}^* \hat{\nabla}^{-\nicefrac{1}{2}}\BRA{\xi}\beta\bigr\rangle
    = \bigl\langle \hat{\nabla}^{-\nicefrac{1}{2}}\BRA{\xi}\beta,
    \hat{\nabla}^{-\nicefrac{1}{2}}\BRA{\xi}\alpha\bigr\rangle
  \end{multline*}
  This implies \(\norm{\xi}^2 K^{-1} =
  (\hat{\nabla}^{-\nicefrac{1}{2}}\BRA{\xi})^*
  \hat{\nabla}^{-\nicefrac{1}{2}}\BRA{\xi} = \KET{\xi}
  \hat{\nabla}^{-1} \BRA{\xi}\).  The operator
  \(K^{-\nicefrac{1}{2}}\) is defined independently of~\(\xi\).
  Hence~\(K\) cannot depend on~\(\xi\).
\end{proof}

\begin{remark}
  \label{rem:unimodular}
  The operator~\(\hat{\nabla}\) generates the modular
  automorphism group of~\(\hat{\QG}\).  Thus~\(\hat{\QG}\) is
  unimodular if and only \(\hat{\nabla}=1\).  In this case, we
  get a scalar formal dimension \(K =
  \norm{\xi}^2\cdot\ID_\Hils\), and \(\Hils_\si=\Hils\) for any
  \emph{irreducible} square-integrable representation.
\end{remark}

\begin{remark}
  \label{rem:compact_orthogonality}
  If the quantum group~\(\QG\) is compact,
  then~\(\hat{\Lambda}\) is the inverse of a bounded map
  \(\Hilsr\to L^1(\QG)\) by
  Remark~\ref{rem:discrete_qg_lambda}.  Therefore, any bounded
  linear functional on~\(\hat{\QG}_\uni\) extends to one
  on~\(\Hilsr\), so that any corepresentation of~\(\QG\) is
  square-integrable.  Thus
  Proposition~\ref{pro:first_orthogonality} and
  Theorem~\ref{the:second_orthogonality} always apply.  The two
  orthogonality relations are part of Woronowicz's Peter--Weyl
  Theorem for compact quantum groups in
  \cites{Woronowicz:Compact_pseudogroups, Woronowicz:CQG}.
\end{remark}

Now we examine the subspace of~\(\Hilsr\) spanned by the
coefficients of an irreducible corepresentation~\(U\)
on~\(\Hils\).

Let~\(\Hils'\) be the Hilbert space completion of \(\dom
K^{-\nicefrac{1}{2}}\subseteq\Hils\) with the conjugate linear
\(\C\)\nb-vector space structure and the inner product
\[
(\xi,\eta) \defeq s(\xi,\eta)
= \langle K^{-\nicefrac{1}{2}}\xi,
K^{-\nicefrac{1}{2}}\eta\rangle.
\]
The Second Orthogonality Relation in
Theorem~\ref{the:second_orthogonality} yields a linear isometry
\[
\Phi\colon \Hils\otimes\Hils'\to\Hilsr,
\qquad \xi\otimes\eta\mapsto \tilde{c}_{\eta,\xi}
= \BRA{\eta}\xi.
\]
We claim that~\(\Phi\) is a bimodule homomorphism with respect
to canonical \(\hat{\QG}\)\nb-bimodule structures on
\(\Hils\otimes\Hils'\) and~\(\Hilsr\).

The bimodule structure on~\(\Hilsr\) is given by the left
regular representation and the right regular
anti-representation:
\begin{equation}
  \label{eq:biregular_rep}
  x\cdot \xi\cdot y \defeq
  \lambda(x)\circ\rho(y)(\xi)
  \qquad\text{for \(x,y\in\hat{\QG}\),
    \(\xi\in\Hilsr\),}
\end{equation}
with~\(\rho\) defined in~\eqref{eq:right_regular_anti}.  The
bimodule structure on \(\Hils\otimes\Hils'\) is defined by
\[
x\cdot (\xi\otimes\eta)\cdot y
\defeq \rho_U(x)(\xi)\otimes
\rho_U\bigl(\hat{\modauto}_{\nicefrac{\ima}{2}}(y)^*\bigr)(\eta).
\]
for \(x\in\hat{\QG}\), \(\xi\in\Hils\), \(\eta\in\Hils'\),
\(y\in\dom\hat{\modauto}_{\nicefrac{\ima}{2}}\).

The equivariance of \(\BRA{\eta}\) established in Lemmas
\ref{lem:KET_equivariant}
and~\ref{lem:act_on_coefficient_right} shows that~\(\Phi\) is
compatible with these bimodule structures.  As a consequence,
the right action of \(\dom\hat{\modauto}_{\nicefrac{\ima}{2}}\)
extends to a \Star{}anti-representation of~\(\hat{\QG}\)
because this happens on~\(\Hilsr\) and~\(\Phi\) is isometric.

We may turn the \(\hat{\QG}\)\nb-bimodule structures above into
\(\QG\)\nb-bicomodule structures, and then~\(\Phi\) is a
homomorphism of \(\QG\)\nb-bicomodules.  More precisely, the
left \(\hat{\QG}\)\nb-module structure corresponds to a
coaction of~\(\QG\), the right \(\hat{\QG}\)\nb-module
structure corresponds to a \emph{left} coaction of~\(\QG\) or,
equivalently, a coaction of the opposite quantum
group~\(\QG^\opp\), that is, of~\(\QG\) with the opposite
comultiplication.  The left and right coactions of~\(\QG\)
commute because the corresponding actions of~\(\hat{\QG}\)
commute.

Let \(P_U\colon \Hilsr\to\Hilsr\) be the projection onto the
range of~\(\Phi\).  Since~\(\Phi\) is a bimodule homomorphism,
its range is a sub-bimodule, so that~\(P_U\) is a bimodule map
with respect to the left and right regular representations.
Equivalently, \(P_U\) belongs to the commutants of both
\(\lambda(\hat{\QG})\) and of \(\rho(\hat{\QG})\).  Since
\(\lambda(\hat{\QG})''\) and \(\rho(\hat{\QG})''\) are commutants
of one another, this means that~\(P_U\) is a central idempotent
in the von Neumann algebra completion \(\lambda(\hat{\QG})''\)
of~\(\lambda(\hat{\QG})\).

It is not hard to see that~\(P_U\) acts on a square-integrable
corepresentation of~\(\QG\) as the projection onto the
\(U\)\nb-isotypical subspace.

Since the representations of~\(\hat{\QG}\) on \(\Hils\)
and~\(\Hils'\) that we use above are both irreducible, the
commutant of the left module structure on
\(\Hils\otimes\Hils'\) is \(\Bound(\Hils')\) acting by
\(T\mapsto 1\otimes T\), and the commutant of the right module
structure is \(\Bound(\Hils)\) acting by \(T\mapsto T\otimes
1\).  Using~\(\Phi\), we see that
\begin{equation}
  \label{eq:range_central_projection}
  P_U \lambda(\hat{\QG})'' \cong \Bound(\Hils),
  \qquad
  P_U \rho(\hat{\QG})'' \cong \Bound(\Hils').
\end{equation}

The First Orthogonality Relation shows that these copies of
\(\Bound(\Hils)\) for non-equivalent irreducible
square-integrable representations are orthogonal.

Conversely, let~\(P\) be a central projection in
\(\lambda(\hat{\QG})''\) with \(P\lambda(\hat{\QG})'' \cong
\Bound(\Hils)\) for some Hilbert space~\(\Hils\), then
\(P\rho(\hat{\QG})''\) is of type~I as well, and the range of a
minimal projection in \(P\rho(\hat{\QG})''\) is an irreducible
representation of \(\lambda(\hat{\QG})''\) because its
commutant is~\(\C\) by construction.  The central projection
associated to the corresponding irreducible square-integrable
corepresentation of~\(\QG\) is exactly~\(P\).  As a result, we
get a bijection between the set of equivalence classes of
irreducible square-integrable corepresentations of~\(\QG\) and
the set of direct summands of type~I in the direct integral
decomposition of the von Neumann algebra
\(\lambda(\hat{\QG})''\).

\smallskip

So far, we have worked rather consistently with representations
of~\(\hat{\QG}_\uni\).  In later sections, we will be equally
consistent in using corepresentations of~\(\QG\).  To see the
connection between the two approaches, we reformulate
Definition~\ref{def:si_Hilbert} in terms of the coaction
\(\coact_U\colon \Hils\to \Mult(\Hils\otimes\QG) \defeq
\Bound(\QG,\Hils\otimes \QG)\):

\begin{lemma}
  \label{lem:coefficient_via_coaction}
  A vector \(\alpha\in\Hils\) is square-integrable if and only
  if \(\langle \coact_U(\alpha),\beta\otimes1_\QG\rangle_\QG\)
  belongs to \(\dom\Lambda\) for all \(\beta\in\Hils\), and in
  this case \(\tilde{c}_{\alpha\beta} = \Lambda\bigl( \langle
  \coact_U(\alpha),\beta\otimes1_\QG\rangle_\QG\bigr)\).
\end{lemma}

Here we write \(\langle x,y\rangle_\QG\defeq x^*\circ y \in
\Bound(\QG,\QG) = \Mult(\QG)\) for \(x,y\in
\Mult(\Hils\otimes\QG)\).

\begin{proof}
  Recall that the map~\(\hat{\Lambda}\) is defined so that
  \(\hat{\Lambda}(x\Haar) = \Lambda(x)\) for all \(x\in\QG\)
  with \(x\Haar\in L^1(\QG)\).  Hence
  \[
  \bigl\langle \tilde{c}_{\alpha\beta},\Lambda(x) \bigr\rangle
  = \bigl\langle \tilde{c}_{\alpha\beta},
  \hat{\Lambda}(x\Haar)\bigr\rangle
  = c_{\alpha,\beta}(x\Haar)
  = \bigl\langle \beta, \rho_U(x\Haar)\alpha \bigr\rangle.
  \]
  By definition,
  \begin{multline*}
    \langle \beta,\rho_U(x\Haar)\alpha\rangle
    = \bigl\langle \beta, (\ID\otimes x\Haar)\coact_U(\alpha)
    \bigr\rangle
    \\= \Haar\bigl(\langle \beta\otimes1_\QG,\coact_U(\alpha)
    \rangle_\QG\cdot x\bigr)
    = \Haar\bigl(\langle\coact_U(\alpha),
    \beta\otimes1_\QG\rangle_\QG^* \cdot x\bigr).
  \end{multline*}
  If \(\langle\coact_U(\alpha), \beta\otimes1_\QG\rangle_\QG
  \in \dom \Lambda\), then we can rewrite this as
  \[
  \langle \beta,\rho_U(x\Haar)\alpha\rangle
  = \bigl\langle \Lambda\bigl(\langle
  \coact_U(\alpha),\beta\otimes1_\QG\rangle_\QG\bigr),
  \Lambda(x)\bigr\rangle,
  \]
  so that~\(\alpha\) is square-integrable and
  \(\tilde{c}_{\alpha\beta} = \Lambda\bigl(\langle
  \coact_U(\alpha),\beta\otimes1_\QG\rangle_\QG\bigr)\).

  Conversely, suppose that~\(\alpha\) is square-integrable and
  let \(y \defeq \langle
  \coact_U(\alpha),\beta\otimes1_\QG\rangle_\QG\) for
  some~\(\beta\).  By assumption, \(x\mapsto \langle
  \beta,\rho_U(x\Haar)\alpha\rangle = \Haar(y^*x)\) is
  bounded with respect to \(\norm{\Lambda(x)}\).  It remains to
  check that this implies \(y\in\dom\Lambda\).

  If \(a\in\dom\Lambda\cap \dom \modauto_{-\ima}\), then
  \(ya\in\dom\Lambda\) as well and
  \[
  \norm{\Lambda(ya)}^2 = \Haar(a^*y^*ya)
  = \Haar\bigl(y^*ya \modauto_{-\ima}(a^*)\bigr)
  \le C\,\bigl\lVert \Lambda\bigl(ya\modauto_{-\ima}(a^*)\bigr)
  \bigr\rVert
  \le C\norm{a}\cdot \norm{\Lambda(ya)}
  \]
  for some constant~\(C\), so that \(\norm{\Lambda(ya)}\le
  C\norm{a}\) for all \(a\in\dom\Lambda\).  Since
  \(\dom\Lambda\) is strictly dense in \(\Mult(\QG)\) we may
  choose an approximate unit \((u_i)_{i\in I}\) in
  \(\dom\Lambda\).  Then \((yu_i)\) converges towards~\(y\).
  Our estimate shows that \(\Lambda(yu_i)\) converges as well.
  Since the operator~\(\Lambda\) is closed, this implies
  \(y\in\dom\Lambda\) as asserted.
\end{proof}

Recall that \(\coact_U(\alpha) = U(\alpha\otimes1_\QG)\).
Hence we can further rewrite
\[
\tilde{c}_{\alpha\beta}
= \Lambda\bigl( \langle
\coact_U(\alpha),\beta\otimes1_\QG\rangle_\QG\bigr)
= \Lambda\bigl( \langle
U(\alpha\otimes1_\QG),\beta\otimes1_\QG\rangle_\QG\bigr).
\]

It is shown in~\cite{Kustermans-Vaes:LCQGvN} that
\(\hat{T}^*\bigl(\Lambda(x)\bigr) =
\Lambda\bigl(S(x^*)\bigr)\), where~\(S\) is the antipode
of~\(\QG\).  Hence the reflection
formula~\eqref{eq:flip_coefficient} is equivalent to
\[
\langle U(\alpha\otimes1),\beta\otimes 1\rangle_\QG
= S\bigl(\langle \alpha\otimes1,
U(\beta\otimes 1)\rangle_\QG\bigr)
\qquad\text{for all \(\alpha,\beta\in\Hils_\si\).}
\]

\section{Integrable coactions}
\label{sec:integrable}

In this section, we define and study integrable coactions of
locally compact quantum groups on \(C^*\)\nb-algebras.  The
main object of interest here is the unbounded map
\((\ID\otimes\Haar)\circ\coact_\coef\), where~\(\Haar\) is the
Haar weight on~\(\QG\) and \(\coact_\coef\in
\Mor(\coef,\coef\otimes\QG)\) is a (right) coaction of a
locally compact quantum group \((\QG,\com)\) on some
\(C^*\)\nb-algebra~\(\coef\), that is, a non-degenerate \Star
homomorphism
\(\coact_\coef\colon\coef\to\Mult(\coef\otimes\QG)\) with
\[
(\coact_\coef\otimes\ID)\circ\coact_\coef
= (\ID\otimes\com)\circ\coact_\coef.
\]

For a group action, \((\ID\otimes\Haar)\circ\coact_\coef\)
becomes the map \(\coef\ni x\mapsto \int_G g\cdot x\,\dd g\),
where the integral is interpreted unconditionally in the strict
topology (see~\cite{Exel:Unconditional}).  Therefore, we call
this map the \emph{averaging map}.  Unless the quantum group
\((\QG,\com)\) is compact, the averaging map is unbounded, and
its domain of definition may be small.  By definition,
integrability means that it is densely defined.

We are going to develop tools to check whether a coaction is
integrable.  They show, in particular, that dual coactions and
stable coactions are integrable.  But first, we rigorously
define the averaging map, using slices with weights as
described in \cite{Kustermans-Vaes:Weight}*{\S3}.  We only
recall this construction briefly here.

The Haar weight is lower semi-continuous by definition.  As
such, it is a limit of an increasing net of bounded weights.
More precisely, we define a partial order for weights
on~\(\QG\) by \(\omega\ll\psi\) if there is \(r\in(0,1)\) with
\(\omega\le r\cdot\psi\).  Then the set
\[
\Sub(\Haar) \defeq \{\psi\in\QG^\dual_+\mid \psi\ll\Haar\}
\]
with partial order~\(\ll\) is a directed set, and
\begin{equation}
  \label{eq:Haar_as_sup}
  \Haar(x) = \lim_{\psi\in\Sub(\Haar)} \psi(x)
\end{equation}
for all \(x\in\QG\).  Even more, \eqref{eq:Haar_as_sup} is used
as a definition to extend~\(\Haar\) to the multiplier algebra
\(\Mult(\QG)\) or to suitable von Neumann algebras.

Given another \(C^*\)\nb-algebra~\(\coef\), we let \(\dom
(\ID\otimes\Haar)^+\) be the set of all \(x\in
\Mult(\coef\otimes\QG)^+\) for which the net
\((\ID\otimes\psi)(x)_{\psi\in\Sub(\Haar)}\) converges in
\(\Mult(\QG)\) in the \emph{strict} topology, and we define
\begin{equation}
  \label{eq:Haar_slice}
  (\ID\otimes\Haar)(x) \defeq
  \lim_{\psi\in\Sub(\Haar)} (\ID\otimes\psi)(x)
  \qquad
  \text{for all \(x\in\dom(\ID\otimes\Haar)^+\).}
\end{equation}
As usual, we let \(\dom(\ID\otimes\Haar)\) be the linear span
of \(\dom(\ID\otimes\Haar)^+\) and extend \(\ID\otimes\Haar\)
linearly to \(\dom(\ID\otimes\Haar)\).
Equation~\eqref{eq:Haar_slice} remains true for \(x\in
\dom(\ID\otimes\Haar)\); but if~\(x\) is not positive, then it
is unclear whether convergence in~\eqref{eq:Haar_slice}
suffices for \(x\in \dom(\ID\otimes\Haar)\).

\begin{definition}
  \label{def:integrable}
  Let \(\coact_\coef\colon \coef\to \Mult(\coef\otimes\QG)\) be
  a coaction of~\(\QG\) on~\(\coef\). We call a positive
  element \(x\in\Mult(\coef)^+\) \emph{integrable} if
  \(\coact_\coef(x)\) belongs to \(\dom (\ID\otimes\Haar)\),
  and we denote the subset of integrable elements in
  \(\Mult(\coef)^+\) by \(\Mult(\coef)_\ible^+\).  Elements in
  \(\Mult(\coef)_\ible \defeq \spann\Mult(\coef)_\ible^+\) are
  also called \emph{integrable}.  We define
  \[
  \Av\colon \Mult(\coef)_\ible \to \Mult(\coef),\qquad
  x \mapsto (\ID\otimes\Haar)\circ\coact_\coef(x),
  \]
  and call this map the \emph{averaging map} for the
  coaction~\(\coact_\coef\) on~\(\coef\).

  We also let \(\coef_\ible^+ \defeq \Mult(\coef)_\ible^+\cap\coef\)
  and \(\coef_\ible \defeq \spann\coef_\ible^+\).
\end{definition}

Although we are mainly interested in the restriction of~\(\Av\)
to~\(\coef\), working with multipliers all the time creates no
additional problems and is sometimes useful.  Even for elements
of~\(\coef\) we use strict convergence to characterise
integrability.  The relevant nets are almost never norm
convergent.

\begin{example}
  \label{exa:free_integrable}
  We consider the standard coaction of a locally compact
  quantum group \((\QG,\com)\) on itself.  More generally, the
  same argument works for \(\coef= D\otimes\QG\) and
  \(\coact_\coef=\ID_D\otimes\com\).  The strong form of the left
  invariance of the Haar weight in
  \cite{Kustermans-Vaes:LCQGvN}*{Proposition 3.1} asserts that
  \[
  (\ID_{D\otimes\QG}\otimes\Haar)\circ\coact_\coef(x) =
  (\ID_D\otimes\Haar)(x)\otimes 1_\QG
  \]
  holds for all \(x\in \Mult(D\otimes\QG)^+\), where both sides
  belong to the extended positive parts of suitable von Neumann
  algebras.  For our purposes, the important point is that one
  side is bounded if and only if the other one is.  Therefore,
  \begin{equation}
    \label{eq:ible_free}
    \Mult(D\otimes\QG)^+_\ible = \dom (\ID\otimes\Haar)^+,\qquad
    \Mult(D\otimes\QG)_\ible = \dom (\ID\otimes\Haar),
  \end{equation}
  and \(\Av(x) = (\ID_D\otimes\Haar)(x)\otimes 1_\QG\) for all
  \(x\in (D\otimes\QG)_\ible\).  In the special case \(D=\C\),
  we get
  \begin{equation}
    \label{eq:ible_free_simple}
    \Mult(\QG)^+_\ible = \dom \Haar^+,\qquad
    \Mult(\QG)_\ible = \dom \Haar,
  \end{equation}
  and \(\Av(x) = 1_\QG\cdot \Haar(x)\) for all \(x\in
  \dom\Haar\).
\end{example}

The properties of the slice map \(\ID\otimes\Haar\) and its
domain listed in~\cite{Kustermans-Vaes:Weight}*{\S3}
immediately imply the assertions in the following lemma:

\begin{lemma}
  \label{lem:integrable}
  The subset of integrable elements \(\Mult(\coef)_\ible\) is a
  hereditary \Star{}subalgebra of \(\Mult(\coef)\), that is, it
  is a \Star{}subalgebra and \(0\le x\le y\) and
  \(y\in\Mult(\coef)_\ible\) imply \(x\in\Mult(\coef)_\ible\).
  If \(x\in\Mult(\coef)_\ible\), then the net
  \((\ID\otimes\psi)\coact_\coef(x)\) for \(\psi\in\Sub(\Haar)\)
  converges strictly towards \(\Av(x)\).  Conversely, if
  \(x\in\Mult(\coef)^+\) is a positive multiplier and the net
  \((\ID\otimes\psi)\coact_\coef(x)\) for \(\psi\in\Sub(\Haar)\)
  converges strictly, then~\(x\) is integrable.
\end{lemma}

\begin{proof}
  Since \(\dom (\ID\otimes\Haar)\) is a hereditary
  \Star{}subalgebra and~\(\coact_\coef\) is a homomorphism,
  \(\Mult(\coef)_\ible\) is a hereditary \Star{}subalgebra.
  The convergence of \((\ID\otimes\psi)\coact_\coef(x)\) towards
  \(\Av(x)\) for \(x\in\Mult(\coef)_\ible\) follows immediately
  from the definitions, see also
  \cite{Kustermans-Vaes:Weight}*{Lemma~3.8}.  The converse
  statement that convergence of this net implies integrability
  holds because a positive element of
  \(\Mult(\coef\otimes\QG)\) belongs to
  \(\dom(\ID\otimes\Haar)\) if and only if it belongs to
  \(\dom(\ID\otimes\Haar)^+\); this uses that
  \(\dom(\ID\otimes\Haar)^+\) is a hereditary cone.
\end{proof}

Recall that a right coaction of~\(\QG\) yields a left
\(\QG^\dual\)\nb-module structure on \(\Mult(\coef)\) via
\(\psi\star x = (\ID\otimes\psi)\coact_ \coef(x)\).  In this
notation, \(\Av(x)\) becomes the limit of \(\psi\star x\) for
\(\psi\in\Sub(\Haar)\), so that we may also write
\(\Av(x)=\Haar\star x\).

Let \(x\in\Mult(\coef)^+\).  Then \((\psi\star
x)_{\psi\in\Sub(\Haar)}\) is an increasing net of positive
multipliers.  For such nets, there are several equivalent ways
to express strict convergence.

First, \((\psi\star x)_{\psi\in\Sub(\Haar)}\) converges
strictly towards~\(y\) if and only if \(b^*\cdot (\psi\star
x)\cdot b\) converges in norm towards \(b^*yb\) for each
\(b\in\coef\) (see \cite{Kustermans-Vaes:Weight}*{Result 3.4}).

Secondly, \((\psi\star x)_{\psi\in\Sub(\Haar)}\) converges
strictly towards \(y\in\Mult(\coef)\) if and only if it
converges weakly, that is, \(\lim \theta(\psi\star x) =
\theta(y)\) for each \(\theta\in \coef^\dual\) (see
\cite{Kustermans-Vaes:Weight}*{Proposition 3.14}).  But mere
weak convergence of the net \(\psi\star x\) in some ambient von
Neumann algebra is not enough: we need the weak limit to be a
multiplier of~\(\coef\) as well.

\begin{example}
  \label{exa:integrable_group}
  Let \(\QG=\CONT_0(G)\) for a locally compact group, with the
  usual comultiplication. Then continuous \(\QG\)\nb-coactions are the
  same as continuous group actions of~\(G\). Bounded positive linear
  functionals on~\(\QG\) correspond to bounded measures
  on~\(G\).  The functionals in \(\Sub(\Haar)\) are those that
  have the form \(f\mapsto \int_G f(g)\cdot w(g)\,\dd{g}\) for a
  function \(w\in L^1(G)\) with \(0\le w\le 1-\varepsilon\) for
  some \(0<\varepsilon<1\); we denote this functional by \(w\,\dd
  g\).  We have \(w\,\dd g\ll w'\,\dd g\) if and only if \(w\le
  (1-\varepsilon)w'\) for some \(0<\varepsilon<1\).

  If~\(\coef\) carries a \(G\)\nb-action and \(x\in\coef\),
  then \(w\,\dd g\star x = \int_G (g\cdot x) w(g)\,\dd g\).
  Integrability means that this net converges strictly; its
  limit is \(\Av(x) = \int_G g\cdot x\,\dd g\).  From this, it
  is easy to see that our definition of integrability agrees
  with the usual one in the group case
  (see~\cite{Rieffel:Integrable_proper}).
\end{example}

\begin{example}
  \label{exa:compact_integrable}
  Let \((\QG,\com)\) be a compact quantum group.  Then the Haar
  weight~\(\Haar\) is itself bounded.  Hence any coaction
  of~\(\QG\) is integrable with \(\coef_\ible=\coef\).
\end{example}

Not surprisingly, the range of the averaging map contains only
\(\QG\)\nb-equivariant multipliers of~\(\coef\):

\begin{lemma}
  \label{lem:av_invariant}
  Let \(x\in\Mult(\coef)_\ible\).  Then
  \(\Av(x)\in\Mult(\coef)\) is \(\QG\)\nb-invariant, that is,
  \(\coact_\coef\bigl(\Av(x)\bigr) = \Av(x)\otimes 1_\QG\).
\end{lemma}

\begin{proof}
  We compute
  \begin{multline*}
    \coact_\coef\bigl(\Av(x)\bigr)
    = (\coact_\coef\otimes\Haar)\coact_\coef(x)
    = (\ID_\coef\otimes\ID_\QG\otimes\Haar)
    (\coact_\coef\otimes\ID_\QG)\coact_\coef(x)
    \\= (\ID_\coef\otimes\ID_\QG\otimes\Haar)
    (\ID_\coef\otimes\com)\coact_\coef(x)
    = (\ID_\coef\otimes\Haar)\coact_\coef(x) \otimes 1_\QG
    = \Av(x) \otimes 1_\QG,
  \end{multline*}
  using the coassociativity of the coaction on~\(\coef\) and
  the left invariance of the Haar measure in the strong form
  \((\ID_\coef\otimes\ID_\QG\otimes\Haar)\circ(\ID_\coef\otimes\com)(x)
  = (\ID_\coef\otimes\Haar)(x)\otimes 1_\QG\) (see
  \cite{Kustermans-Vaes:LCQGvN}*{Proposition 3.1}).  To justify
  the first three equalities, we should replace~\(\Haar\) by a
  limit over \(\psi\in\Sub(\Haar)\).
\end{proof}

Next we study how integrability behaves with respect to the
action of~\(\QG^\dual\).  For group actions, the subspace of
integrable elements is \(G\)\nb-invariant; but the
\(G\)\nb-action on the space of integrable elements does not
integrate to a module structure over \(L^1(G)\) because the
modular function intervenes.  Therefore, it is to be expected
that the modular element of the locally compact quantum group
\((\QG,\com)\) appears at this point.  Recall that the
\emph{modular element}~\(\modular\) is an unbounded multiplier
affiliated with~\(\QG\) that is characterised by the property
that the unbounded weight
\[
\Haar_\modular =
\modular^{\nicefrac{1}{2}}\Haar\modular^{\nicefrac{1}{2}}\colon
x\mapsto \Haar(\modular^{\nicefrac{1}{2}} x
\modular^{\nicefrac{1}{2}})
\]
is a right invariant Haar weight on~\(\QG\), that is,
\(\Haar_\modular \star \omega = \Haar_\modular\cdot \omega(1)\)
for all \(\omega\in\QG^\dual\).  The modular element is
group-like, that is, \(\com(\modular) =
\modular\otimes\modular\).  This implies \(\com(\modular^t) =
\modular^t\otimes\modular^t\) for all \(t\in\R\).

\begin{lemma}
  \label{lem:integrable_invariant}
  Let \(x\in\Mult(\coef)_\ible\) and let \(\omega\in\QG^\dual\)
  be such that the weight \(\omega_\modular =
  \modular^{\nicefrac{1}{2}}\omega\modular^{\nicefrac{1}{2}}\)
  is bounded as well.  Then \(\omega\star x\in
  \Mult(\coef)_\ible\) and \(\Av(\omega\star x) =
  \omega(\modular)\cdot \Av(x)\).
\end{lemma}

\begin{proof}
  We may assume without loss of generality that \(x\)
  and~\(\omega\) are positive, the general case follows by
  polarisation.  Then \(\omega\star x\) is also positive.  Thus
  we must show that the net \(\psi\star(\omega\star x)\) for
  \(\psi\in\Sub(\Haar)\) converges towards
  \(\omega(\modular)\Av(x) = \omega(\modular)\cdot (\Haar\star
  x)\).  Since the convolution~\(\star\) is associative, this
  boils down to the convergence of weights
  \[
  \lim_{\psi\in\Sub(\Haar)} \psi\star\omega =
  \omega(\modular)\cdot\Haar.
  \]
  Let \(y\in\QG^+\) satisfy
  \(\modular^{\nicefrac{1}{2}}y\modular^{\nicefrac{1}{2}}\in
  \dom\Haar \subseteq \QG\).  Using
  \(\com(\modular^{\nicefrac{1}{2}}) =
  \modular^{\nicefrac{1}{2}}\otimes
  \modular^{\nicefrac{1}{2}}\) and the right invariance
  of~\(\Haar_\modular\), we get
  \begin{multline*}
    \lim_{\psi\in\Sub(\Haar)}
    \psi\star\omega
    (\modular^{\nicefrac{1}{2}}y\modular^{\nicefrac{1}{2}})
    = \lim_{\psi\in\Sub(\Haar)}
    (\psi_\modular\otimes\omega_\modular)\bigl(\com(y)\bigr)
    \\= \omega_\modular\circ(\Haar_\modular\otimes\ID)\com(y)
    = \omega_\modular\bigl(1_\QG\cdot\Haar_\modular(y)\bigr)
    = \omega(\modular)\cdot
    \Haar(\modular^{\nicefrac{1}{2}}y\modular^{\nicefrac{1}{2}}).
  \end{multline*}

  Hence \((\psi\star\omega)_{\psi\in\Sub(\Haar)}\) converges
  towards \(\omega(\modular)\cdot\Haar\) on a dense subset
  of~\(\QG\).  This implies the assertion because this net of
  weights is increasing.
\end{proof}

\begin{lemma}
  \label{lem:integrable_coaction}
  Let~\(\coact_\coef\) be a coaction of a locally compact quantum
  group~\(\QG\) on a \(C^*\)\nb-algebra~\(\coef\).  The
  following statements are equivalent:
  \begin{enumerate}[label=\textup{(\arabic{*})}]
  \item \(\coef_\ible\) is norm dense in~\(\coef\);

  \item \(\coef^+_\ible\) is norm dense in~\(\coef^+\);

  \item \(\Mult(\coef)_\ible\) is strictly dense
    in~\(\Mult(\coef)\);

  \item \(\Mult(\coef)^+_\ible\) is strictly dense
    in~\(\Mult(\coef)^+\);
  \end{enumerate}

  If, in addition, \(\coef\) is \(\sigma\)\nb-unital, then
  \textup{(1)}--\textup{(4)} are also equivalent to:
  \begin{enumerate}[label=\textup{(\arabic{*})}]
  \item[\textup{(5)}] \(\coef^+_\ible\) contains a strictly
    positive element of~\(\coef\);

  \item[\textup{(6)}] there exists a strictly positive
    integrable multiplier.
  \end{enumerate}
\end{lemma}

\begin{proof}
  The equivalences (1)\(\iff\)(2) and (3)\(\iff\)(4) follow
  because \(\Mult(\coef)^+_\ible\) is a hereditary cone in
  \(\Mult(\coef)\), and (2)\(\Longrightarrow\)(4) follows
  because~\(\coef\) is strictly dense in \(\Mult(\coef)\) and
  the strict topology is weaker than the norm topology.

  Now we check that, conversely, (4)\(\Longrightarrow\)(2).
  Let \(x\in\Mult(\coef)^+_\ible\) and \(y\in\coef^+\), then
  \(x^{\nicefrac{1}{2}}yx^{\nicefrac{1}{2}} \le
  x^{\nicefrac{1}{2}}\norm{y}x^{\nicefrac{1}{2}} =
  \norm{y}\cdot x\), so that
  \(x^{\nicefrac{1}{2}}yx^{\nicefrac{1}{2}}\) is integrable as
  well.  If \(\Mult(\coef)^+_\ible\) is strictly dense in
  \(\Mult(\coef)^+\), then the set
  \(\{x^{\nicefrac{1}{2}}yx^{\nicefrac{1}{2}} \mid x\in
  \Mult(\coef)^+_\ible, y\in \coef^+\}\) is norm-dense
  in~\(\coef^+\) and consists of integrable elements.  Thus~(4)
  implies~(2).

  A similar argument shows that~(6) implies~(2) because
  \(x^{\nicefrac{1}{2}} \coef x^{\nicefrac{1}{2}}\) is dense
  in~\(\coef\) if~\(x\) is a strictly positive multiplier
  of~\(\coef\).  Since (5)\(\Longrightarrow\)(6) is trivial, it
  remains to check that~(2) implies~(5) if~\(\coef\) is
  \(\sigma\)\nb-unital.  This ensures that~\(B\) contains a
  strictly positive element~\(h\).  By integrability, there is
  a sequence \((x_n)_{n\in\N}\) of positive integrable elements
  that converges towards~\(h\).  Then \(\sum \lambda_n x_n\) is
  strictly positive for any sequence of strictly positive
  numbers \((\lambda_n)\).  Choose~\((\lambda_n)\) to decay so
  rapidly that both \(\sum \lambda_n \norm{x_n}\) and \(\sum
  \lambda_n \norm{\Av(x_n)}\) remain bounded.  Then \(\sum
  \lambda_n x_n\) is again integrable, and it is also strictly
  positive.
\end{proof}

\begin{definition}
  \label{def:integrable_coaction}
  The coaction~\(\coact_\coef\) on~\(\coef\) is called
  \emph{integrable} if the equivalent conditions in
  Lemma~\ref{lem:integrable_coaction} are satisfied.
\end{definition}

\begin{proposition}
  \label{pro:integrable_functorial}
  Let \(\coef\) and~\(\coefd\) be \(C^*\)\nb-algebras with
  coactions \(\coact_\coef\) and~\(\coact_\coefd\) of
  \((\QG,\com)\), and let \(\pi\colon
  \Mult(\coef)\to\Mult(\coefd)\) be a \(\QG\)\nb-equivariant,
  strictly continuous, positive linear map that is
  non-degenerate in the sense that \(\pi(\coef^+)\cdot \coefd\)
  is dense in~\(\coefd\).

  Then \(\pi\bigl(\Mult(\coef)_\ible\bigr) \subseteq
  \Mult(\coefd)_\ible\), and \(\pi\bigl(\Av(x)\bigr) =
  \Av\bigl(\pi(x)\bigr)\) for all \(x\in\Mult(\coef)_\ible\).
  As a consequence, if~\(\coef\) is integrable, so
  is~\(\coefd\).
\end{proposition}

In most applications, \(\pi\) is the strictly continuous
extension of a non-degenerate equivariant \Star{}homomorphism
\(\coef\to\Mult(\coefd)\).  We allow more general maps to
prepare for Proposition~\ref{pro:integrable_via_cutoff}, which
characterises integrable coactions by the existence of an
equivariant, non-degenerate, completely positive map
\(\QG\to\Mult(\coef)\).

\begin{proof}
  Since~\(\pi\) is \(\QG\)\nb-equivariant, \(\pi(\psi\star
  x)=\psi\star \pi(x)\) for all \(x\in\Mult(\coef)\),
  \(\psi\in\Sub(\Haar)\).  Since~\(\pi\) is strictly continuous
  and positive, this shows that \(\pi(x)\) is integrable
  once~\(x\) is and that~\(\pi\) intertwines the averaging maps
  for \(\coef\) and~\(\coefd\).

  Now let \((\coef,\coact_\coef)\) be integrable, that is,
  \(\coef_\ible^+\) is dense in~\(\coef^+\).  Then the set of
  \(\pi(b)d\pi(b)\) with \(b\in\coef_\ible^+\),
  \(d\in\coefd^+\) is dense in~\(\coefd^+\) because~\(\pi\) is
  non-degenerate.  Since \(\pi(\coef_\ible^+)\subseteq
  \Mult(\coefd)_\ible^+\) and \(\Mult(\coefd)_\ible^+\) is a
  hereditary cone, \(\pi(b)d\pi(b) \in \Mult(\coefd)^+_\ible\)
  for all \(b\in\coef_\ible^+\), \(d\in\coefd^+\).  Hence
  \(\coefd_\ible^+ = \coefd^+\cap \Mult(\coefd)_\ible^+\) is
  dense in~\(\coefd^+\), so that \((\coefd,\coact_\coefd)\) is
  integrable.
\end{proof}

Most of the remaining results in this section are applications
of Proposition~\ref{pro:integrable_functorial}.

\begin{corollary}
  \label{cor:integrable_ideal_quotient}
  Let~\(\coef\) be a \(C^*\)\nb-algebra with a
  coaction~\(\coact_\coef\) of \((\QG,\com)\) and let
  \(I\subseteq \coef\) be a \(\QG\)\nb-invariant ideal.
  Equip~\(I\) and the quotient \(\coef/I\) with the induced
  coactions of~\(\QG\).  If \((\coef,\coact_\coef)\) is
  integrable, so are \((I,\coact_I)\) and
  \((\coef/I,\coact_{\coef/I})\).
\end{corollary}

\begin{proof}
  The \Star{}homomorphisms \(\coef\to\Mult(I)\) and \(\coef\to
  \coef/I\) are non-degenerate and equivariant.  Hence the
  assertion follows from
  Proposition~\ref{pro:integrable_functorial}.
\end{proof}

Conversely, suppose that the coactions on \(I\) and~\(\coef/I\)
are integrable.  If the extension \(I\into \coef\prto \coef/I\)
has an equivariant completely positive section, then
Proposition~\ref{pro:integrable_functorial} applied to the
resulting linear isomorphism \(I\oplus \coef/I\to \coef\) shows
that \((\coef,\coact_\coef)\) is integrable.

But without extra assumption, an extension of two integrable
coactions need not be integrable.  Group actions on commutative
\(C^*\)\nb-algebras provide counterexamples.  It is shown by
Marc Rieffel in~\cite{Rieffel:Integrable_proper} that a group
action on a locally compact space~\(X\) is proper if and only
if the induced group action on \(\CONT_0(X)\) is integrable.
But it is possible to get non-proper actions by glueing
together proper actions.

\begin{corollary}
  \label{cor:corner_integrable}
  Let \((\coef,\coact_\coef)\) be an integrable coaction of
  \((\QG,\com)\) and let \(p\in\Mult(\coef)\) be a
  \(\QG\)\nb-invariant projection, that is,
  \(\coact_\coef(p)=p\otimes 1\). Then the coaction
  on~\(\coef\) restricts to a coaction on the \emph{corner
    subalgebra} \(p\coef p\), which is again integrable.
\end{corollary}

\begin{proof}
  Apply Proposition~\ref{pro:integrable_functorial} to the map
  \(x\mapsto pxp\) from~\(\coef\) to~\(p\coef p\), which is
  completely positive, equivariant, and non-degenerate.
\end{proof}

Corner subalgebras are a special case of \emph{hereditary
  subalgebras}.  But integrability does not appear to be
inherited by hereditary subalgebras in general.

We are going to show that dual coactions are integrable.  First
we briefly recall the definition.  The two articles
\cites{Baaj-Skandalis:Unitaires, Vaes:Induction_Imprimitivity}
use different conventions at this point: the dual quantum group
is defined using the right instead of the left regular
corepresentation in~\cite{Baaj-Skandalis:Unitaires}.  Here we
follow the conventions in \cites{Kustermans-Vaes:LCQG,
  Vaes:Induction_Imprimitivity}, so that we get the
\(C^*\)\nb-commutant~\(\hat{\QG}^\comm\) instead
of~\(\hat{\QG}\), which appears
in~\cite{Baaj-Skandalis:Unitaires}.

Let \((\QG,\com)\) be a locally compact quantum group,
let~\(\coef\) be a \(C^*\)\nb-algebra, and let
\(\coact_\coef\colon \coef\to \Mult(\coef\otimes\QG)\) be a
coaction of \((\QG,\com)\) on~\(\coef\).  So far, the only
assumption on~\(\coact_\coef\) that we needed was that it is
non-degenerate and coassociative.  To define a crossed product,
we also require the coaction to be \emph{continuous} in the sense that
\begin{equation}
  \label{eq:coaction_continuous}
  \clspan \bigl( (1_\coef\otimes\QG)\cdot \coact_\coef(\coef)\bigr)
  = \coef\otimes\QG.
\end{equation}

We represent \(\QG\) and~\(\hat{\QG}\) on the Hilbert
space~\(\Hilsr\) as usual, and view \(\coef\otimes\QG\) as a
\(C^*\)\nb-algebra of adjointable operators on the Hilbert
\(\coef\)\nb-module \(\coef\otimes\Hilsr\).  This extends to a
\Star{}representation \(\Mult(\coef\otimes\QG) \to
\Bound(\coef\otimes\Hilsr)\).  Let
\begin{equation}
  \label{eq:Cstar_commutant_QG}
  \hat{\QG}^\comm \defeq \hat{J}\hat{\QG}\hat{J}
  \subseteq\Bound(\Hilsr)\subseteq\Bound(\coef\otimes\Hilsr)
\end{equation}
be the \(C^*\)\nb-commutant of~\(\hat{\QG}\).  The continuity
assumption~\eqref{eq:coaction_continuous} implies that
\begin{equation}
  \label{eq:definition_crossed}
  \coef\rtimes_\red \hat{\QG}^\comm \defeq
  \clspan \bigl( \coact_\coef(\coef)\cdot
  (1\otimes\hat{\QG}^\comm)\bigr)
  = \clspan \bigl((1\otimes\hat{\QG}^\comm)
  \cdot\coact_\coef(\coef)\bigr)
\end{equation}
is a \(C^*\)\nb-subalgebra of \(\Bound(\coef\otimes\Hilsr)\).
This is the \emph{reduced \(C^*\)\nb-crossed product} of the
coaction~\(\coact_\coef\) on~\(\coef\).  The \emph{dual
  coaction}~\(\hat{\coact}_\coef\) of \(\hat{\QG}^\comm\) on
\(\coef\rtimes_\red \hat{\QG}^\comm\) is defined by
\[
\hat{\coact}_\coef\bigl(\coact_\coef(a)\cdot (1\otimes x)\bigr)
\defeq (\coact_\coef(a)\otimes 1_\QG)\cdot \bigl(1_\coef\otimes
\hat{\com}^\comm(x)\bigr),
\]
where~\(\hat{\com}^\comm\) denotes the canonical
comultiplication on~\(\hat{\QG}^\comm\).

Besides the reduced crossed product, there is the \emph{full
  \(C^*\)\nb-crossed product} \(\coef\rtimes \hat{\QG}^\comm_\uni\),
which is defined to be universal for covariant representations
of \((\coef,\coact_\coef)\) and \((\QG,\com)\).  The two
crossed products are linked to their constituents by canonical
morphisms
\begin{alignat*}{2}
  j_\coef&\in\Mor(\coef,\coef\rtimes \hat{\QG}^\comm_\uni),&\qquad
  j_{\hat{\QG}^\comm_\uni}&\in\Mor(\hat{\QG}^\comm_\uni,
  \coef\rtimes \hat{\QG}^\comm_\uni),\\
  j_\coef^\red&\in\Mor(\coef,\coef\rtimes_\red \hat{\QG}^\comm),&\qquad
  j_{\hat{\QG}^\comm}^\red&\in\Mor(\hat{\QG}^\comm,
  \coef\rtimes_\red \hat{\QG}^\comm),
\end{alignat*}
where \(\Mor(A,B)\) denotes the set of non-degenerate
\Star{}homomorphisms \(A\to\Mult(B)\) and
\((\QG_\uni,\com_\uni)\) is the universal locally compact
quantum group associated to the reduced quantum group
\((\QG,\com)\), see~\cite{Kustermans:LCQG_universal}.

\begin{proposition}
  \label{pro:dual_coaction_integrable}
  Let \((\QG,\com)\) be a locally compact quantum group,
  let~\(\coef\) be a \(C^*\)\nb-algebra, and let
  \(\coact_\coef\colon \coef\to \Mult(\coef\otimes
  \hat{\QG}^\comm)\) be a continuous coaction of
  \((\hat{\QG}^\comm,\hat{\com}^\comm)\).  The associated dual
  coactions of \((\QG,\com)\) on the reduced and full
  \(C^*\)\nb-crossed products \(\coef\rtimes_\red\QG\) and
  \(\coef\rtimes\QG_\uni\) are both integrable.
\end{proposition}

\begin{proof}
  The canonical morphisms \(\QG\to \Mult(\coef\rtimes_\red\QG)\) and
  \(\QG_\uni\to \Mult(\coef\rtimes\QG_\uni)\) are
  \(\QG\)\nb-equivariant non-degenerate \Star{}homomorphisms;
  here we equip \(\QG\) and~\(\QG_\uni\) with the coactions
  \(\com\) and \((\ID\otimes\pi)\circ\com_\uni\), where
  \(\pi\colon \QG_\uni\to\QG\) is the canonical map, and the
  crossed products carry the dual coactions.  By
  Proposition~\ref{pro:integrable_functorial}, it therefore
  suffices to check the integrability of the coactions on
  \(\QG\) and~\(\QG_\uni\).  The integrability of
  \((\QG,\com)\) is Example~\ref{exa:free_integrable}.  Similar
  assertions hold for~\(\QG_\uni\) instead of~\(\QG\): it turns
  out that the following assertions are equivalent for
  \(x\in\QG_\uni\):
  \begin{itemize}
  \item \(x\) is integrable;

  \item \(\pi(x) \in \QG\) is integrable;

  \item \(\pi(x) \in \dom\Haar\);

  \item \(x \in \dom(\Haar\circ\pi)\).
  \end{itemize}
  Hence the canonical coaction on~\(\QG_\uni\) is integrable as
  well.
\end{proof}

The next class of examples we consider are \emph{stable
  coactions}.  Let~\(\coact_\coef\) be a right coaction of
\((\QG,\com)\) on a \(C^*\)\nb-algebra~\(\coef\).  We equip the
Hilbert space~\(\Hilsr\) with the left regular corepresentation
of~\(\QG\) defined in~\eqref{eq:coact_lambda}.  Then there is a
tensor product coaction of \((\QG,\com)\) on the Hilbert
\(\coef\)\nb-module \(\coef\otimes\Hilsr\), defined by
\begin{equation}
  \label{eq:standard_regular_coact}
  \coact_{\coef\otimes\Hilsr}(x) =
  (1\otimes\Sigma W)(\coact_\coef\otimes\ID_{\Hilsr})(x)
\end{equation}
for all \(x\in \coef\otimes\Hilsr\), where \(\Sigma\colon
\QG\otimes\Hilsr\to\Hilsr\otimes\QG\) is the coordinate flip
and~\(W\) is the multiplicative unitary for \((\QG,\com)\)
described in~\eqref{eq:mult_unitary}; this tensor product
coaction is a special case of
\cite{Baaj-Skandalis:Hopf_KK}*{Proposition 2.10}.

The coaction on \(\coef\otimes\Hilsr\) induces a coaction on
the \(C^*\)\nb-algebra of compact operators
\(\Comp(\coef\otimes\Hilsr)\cong \coef\otimes \Comp(\Hilsr)\).
Coactions of this form are called \emph{stable}.

\begin{proposition}
  \label{pro:stable_integrable}
  Stable coactions are integrable.
\end{proposition}

\begin{proof}
  We equip~\(\QG\) with the standard coaction~\(\com\) of
  itself as in Example~\ref{exa:free_integrable}.  The standard
  representation \(\QG\to\Bound(\Hilsr)\) is
  \(\QG\)\nb-equivariant with respect to the \emph{right}
  regular corepresentation~\(V\) of~\(\QG\) on~\(\Hilsr\)
  because \(\com(x) = V(x\otimes 1)V^*\) for all \(x\in\QG\)
  (see \cite{Kustermans-Vaes:LCQGvN}*{page 86}).  The left
  regular corepresentation is described, as a right
  corepresentation, by the unitary~\(\hat{W}^*\).  The unitary
  \(U=J \hat{J}\) that combines the modular conjugations
  of~\(\QG\) and~\(\hat{\QG}\) provides a unitary equivalence
  between the left and right regular corepresentations, that
  is, \(V=(U^*\otimes\ID)\hat{W}^*(U\otimes\ID)\) (see
  \cite{Kustermans-Vaes:LCQGvN}*{Proposition 2.15}).

  Hence we get a \(\QG\)\nb-equivariant, non-degenerate
  \Star{}homomorphism
  \[
  \QG\to\Bound(\coef\otimes\Hilsr)
  = \Mult\bigl(\Comp(\coef\otimes\Hilsr)\bigr),
  \qquad x\mapsto \ID_\coef\otimes UxU^*.
  \]
  The coaction~\(\com\) on~\(\QG\) is integrable by
  Example~\ref{exa:free_integrable}.  Finally,
  Proposition~\ref{pro:integrable_functorial} shows that the
  coaction on~\(\Comp(\coef\otimes\Hilsr)\) is integrable.
\end{proof}

We have chosen to work with coactions of reduced locally
compact quantum groups \((\QG,\com)\).  Of course, since its
universal companion \((\QG_\uni,\com_\uni)\) also has Haar
weights, we could also define and work with integrability in
the universal setting.  The following result shows that both
approaches are equivalent.

\begin{proposition}
  Let \(\coact_\coef^\uni\colon
  \coef\to\Mult(\coef\otimes\QG_\uni)\) be a coaction of
  \((\QG_\uni,\com_\uni)\) on a \(C^*\)\nb-algebra~\(\coef\)
  and consider its \emph{reduced form}
  \[
  \coact_\coef \defeq (\ID\otimes\pi)\circ\coact_\coef^\uni
  \colon \coef\to\Mult(\coef\otimes\QG)
  \]
  where \(\pi\colon \QG_\uni\to\QG\) is the canonical map.
  Then \((\coef,\coact_\coef^\uni)\) is integrable if and only
  if \((\coef,\coact_\coef)\) is.  Moreover, the spaces of
  integrable elements in \(\Mult(\coef)\) with respect to both
  coactions coincide and \(\Av(x)=\Av_\uni(x)\) for every
  integrable element \(x\), where \(\Av_\uni\) denotes the
  averaging map for \((\coef,\coact_\coef^\uni)\).
\end{proposition}

\begin{proof}
  Recall that the left invariant Haar weight on~\(\QG_\uni\) is
  given by \(\Haar_\uni(x)=\Haar\bigl(\pi(x)\bigr)\).  It
  follows that \(x\in\dom(\Haar_\uni)\) if and only if
  \(\pi(x)\in \dom(\Haar)\).  This and
  \cite{Kustermans-Vaes:Weight}*{Proposition~3.14} show that
  \(x\in \dom(\ID_\coef\otimes\Haar_\uni)\) if and only if
  \((\ID_\coef\otimes\pi)(x)\in\dom(\ID_\coef\otimes\Haar)\),
  and in this case
  \[
  (\ID_\coef\otimes\Haar_\uni)(x)
  = (\ID_\coef\otimes\Haar)\bigl((\ID_\coef\otimes\pi)(x)\bigr).
  \]
  Now the assertions follow.
\end{proof}

\section{Square-integrable coactions on Hilbert modules}
\label{sec:si_Hilm}

Square-integrable coactions on Hilbert modules contain
square-integrable corepresentations on Hilbert spaces and
integrable coactions on \(C^*\)\nb-algebras as special cases.
Thus they combine the two situations studied in Sections
\ref{sec:si_Hilbert_space} and~\ref{sec:integrable}.

First we define coactions on Hilbert modules, following Saad
Baaj and Georges Skandalis \cite{Baaj-Skandalis:Hopf_KK}*{\S2}.
Then we define square-integrable vectors in Hilbert modules and
associate operators \(\BRA{\xi}\) and \(\KET{\xi}\) to them.
We provide some examples of square-integrable Hilbert modules
and study the general properties of the space of
square-integrable vectors and the operators \(\BRA{\xi}\) and
\(\KET{\xi}\).

\subsection{Quantum group coactions on Hilbert modules}
\label{sec:coactions_Hilm}

Let \((\QG,\com)\) be a locally compact quantum group and
let~\(\coef\) be a \(C^*\)\nb-algebra with a
coaction~\(\coact_\coef\) of~\(\QG\).  Let~\(\Hilm\) be a
Hilbert \(\coef\)\nb-module.  We recall the definition and some
basic properties of coactions on Hilbert modules.  We refer
to~\cite{Baaj-Skandalis:Hopf_KK} for further details.

\begin{definition}
  \label{def:coaction_Hilm}
  A \emph{coaction} of \((\QG,\com)\) on~\(\Hilm\) is a linear
  map
  \[
  \coact_\Hilm\colon \Hilm\to \Mult(\Hilm\otimes\QG) \defeq
  \Bound(\coef\otimes\QG,\Hilm\otimes\QG)
  \]
  that satisfies the following conditions:
  \begin{itemize}
  \item \(\coact_\Hilm(\xi\cdot b) =
    \coact_\Hilm(\xi)\cdot\coact_\coef(b)\) for all
    \(\xi\in\Hilm\), \(b\in\coef\);

  \item \(\bigl\langle \coact_\Hilm(\xi),\coact_\Hilm(\eta)
    \bigr\rangle_{\Mult(\coef\otimes\QG)} \defeq
    \coact_\Hilm(\xi)^*\circ\coact_\Hilm(\eta) =
    \coact_\coef\bigl(\langle\xi,\eta\rangle_\coef\bigr)\) for
    all \(\xi,\eta\in\Hilm\);

  \item \(\clspan \coact_\Hilm(\Hilm)\cdot (\coef\otimes\QG) =
    \Hilm\otimes\QG\);

  \item \((\coact_\Hilm\otimes\ID_\QG)\circ\coact_\Hilm =
    (\ID_\Hilm\otimes\com)\circ\coact_\Hilm\).

  \end{itemize}
  The last condition involves extensions of
  \(\coact_\Hilm\otimes\ID_\QG\) and \(\ID_\Hilm\otimes\com\)
  to the multiplier modules, which exist in the presence of the
  other conditions.

  We also call a Hilbert \(\coef\)\nb-module~\(\Hilm\) with
  such a coaction of \((\QG,\com)\) a
  \emph{\(\QG\)\nb-equivariant Hilbert \(\coef\)\nb-module}.
\end{definition}

\begin{example}
  \label{exa:Cstar_equivariant_Hilm}
  The \(C^*\)\nb-algebra~\(\coef\), viewed as a Hilbert module
  over itself becomes a \(\QG\)\nb-equivariant Hilbert
  \(\coef\)\nb-module for the given coaction~\(\coact_\coef\).
\end{example}

\begin{remark}
  \label{rem:discontinuous}
  We follow \cites{Baaj-Skandalis-Vaes:Non-semi-regular,
    Vaes:Induction_Imprimitivity} and deviate
  from~\cite{Baaj-Skandalis:Hopf_KK} by not requiring
  \(\coact_\Hilm(\Hilm)\cdot (1\otimes\QG)\) and
  \((1\otimes\QG)\cdot\coact_\Hilm(\Hilm)\) to be contained in
  \(\Hilm\otimes\QG\).  Many of the basic assertions on Hilbert
  module coactions in~\cite{Baaj-Skandalis:Hopf_KK} do not
  require this conditions.

  It is necessary to drop this condition in order to treat
  non-regular quantum groups because a locally compact quantum
  group~\(\QG\) is regular if and only if the
  coaction~\(\coact_\lambda\) on~\(\Hilsr\) from the left
  regular corepresentation satisfies the condition
  \((1\otimes\QG)\cdot\coact_\lambda(\Hilsr)\subseteq
  \Hilsr\otimes\QG\).
\end{remark}

A coaction on~\(\Hilm\) induces a coaction on the
\(C^*\)\nb-algebra of compact operators
\[
\coact_{\Comp(\Hilm)}\colon \Comp(\Hilm) \to
\Mult(\Comp(\Hilm)\otimes\QG),
\qquad
\ket{\xi}\bra{\eta} \mapsto
\coact_\Hilm(\xi)\circ\coact_\Hilm(\eta)^*.
\]
Conversely, a coaction of~\(\QG\) on the (linking) \(C^*\)\nb-algebra
\(\Comp(\Hilm\oplus\coef)\) that restricts to the given
coaction on the corner \(\Comp(\coef)\cong\coef\) is equivalent
to a coaction of \((\QG,\com)\) on~\(\Hilm\).

The coaction on \(\Comp(\Hilm)\) extends to a map
\[
\coact_{\Bound(\Hilm)} \colon \Bound(\Hilm) \to
\Bound(\Hilm\otimes\QG)
\]
that satisfies
\begin{equation}
  \label{eq:coact_on_operator}
  \coact_\Hilm(T\xi) =
  \coact_{\Bound(\Hilm)}(T)\coact_\Hilm(\xi)
  \qquad\text{for all \(T\in\Bound(\Hilm)\), \(\xi\in\Hilm\).}
\end{equation}
This extension is constructed easily using a third equivalent
description of a coaction on~\(\Hilm\) in terms of a unitary
operator
\begin{equation}
  \label{eq:coaction_via_unitary}
  V\colon \Hilm\otimes_{\coact_\coef} (\coef\otimes\QG) \to
  \Hilm\otimes\QG.
\end{equation}
We simply put
\begin{equation}
  \label{eq:coaction_on_operators}
  \coact_{\Bound(\Hilm)}(T) \defeq
  V(T\otimes_{\coact_\coef}\ID_{\coef\otimes\QG})V^*
\end{equation}
Similarly, the coaction~\(\coact_\Hilm\) extends to a map
\[
\coact_{\Mult(\Hilm)}\colon \Mult(\Hilm)
\defeq \Bound(\coef,\Hilm)
\to \Mult(\Hilm\otimes\QG) \defeq
\Bound(\coef\otimes\QG,\Hilm\otimes\QG)
\]
via \(T\mapsto V\circ
(T\otimes_{\coact_\coef}\ID_{\coef\otimes\QG})\) and the
identification \(\coef\otimes_{\coact_\coef} (\coef\otimes\QG)
\cong \coef\otimes\QG\).

\subsection{Square-integrable elements of equivariant Hilbert
  modules}
\label{sec:si_Hilm_def}

The following definition generalises the description of
square-integrable vectors for Hilbert space representations in
Lemma~\ref{lem:coefficient_via_coaction}.  It uses the slice
map
\[
\ID\otimes\Lambda\colon
\Mult(\coef\otimes\QG) \supseteq \dom (\ID\otimes\Lambda)\to
\Mult(\coef\otimes\Hilsr)
\]
associated to the map \(\Lambda\colon \QG\supseteq \dom
\Lambda\to\Hilsr\) defined
in~\cite{Kustermans-Vaes:LCQG}*{\S1.5}.  We will describe this
map in more detail below.

\begin{definition}
  \label{def:si_Hilm}
  Let~\(\Hilm\) be a \(\QG\)\nb-equivariant Hilbert
  \(\coef\)\nb-module with coaction~\(\coact_\Hilm\), let
  \(V\colon \Hilm\otimes_{\coact_\coef} (\coef\otimes\QG) \to
  \Hilm\otimes\QG\) be the unitary operator that
  describes~\(\coact_\Hilm\).  We call \(\xi\in\Hilm\) or, more
  generally, \(\xi\in\Mult(\Hilm)\) \emph{square-integrable} if
  \[
  \coact_{\Mult(\Hilm)}(\xi)^* (\eta\otimes1_\QG) =
  (\xi\otimes_{\coact_\coef} \ID_{\coef\otimes\QG})^*V^*
  (\eta\otimes1_\QG)
  \in \Mult(\coef\otimes\QG)
  \]
  belongs to the domain of \(\ID_\coef\otimes\Lambda\) for all
  \(\eta\in\Hilm\), and we define
  \[
  \BRA{\xi}\eta \defeq (\ID\otimes\Lambda)
  \bigl(\coact_{\Mult(\Hilm)}(\xi)^* (\eta\otimes1_\QG)\bigr).
  \]
  We let \(\Mult(\Hilm)_\si\) and~\(\Hilm_\si\) be the
  subspaces of square-integrable elements.  We call
  \((\Hilm,\coact_\Hilm)\) \emph{square-integrable}
  if~\(\Hilm_\si\) is norm-dense in~\(\Hilm\).
\end{definition}

To give meaning to this definition, we must describe the slice
map \(\ID_\coef\otimes\Lambda\) and its domain.  We do this
slightly more directly than in
\cite{Kustermans-Vaes:Weight}*{\S3}.

Recall that the Haar weight~\(\Haar\) is the limit of the
increasing net of bounded weights
\((\psi)_{\psi\in\Sub(\Haar)}\) and that
\begin{multline*}
  \dom (\ID_\coef\otimes\Haar)^+ \defeq
  \bigl\{x\in \Mult(\coef\otimes\QG)^+ \bigm|
  \\\text{the net
    \(\bigl((\ID\otimes\psi)(x)\bigr)_{\psi\in\Sub(\Haar)}\) in
    \(\Mult(\coef)\) converges strictly}\bigr\}.
\end{multline*}
We define
\[
\dom (\ID_\coef\otimes\Lambda) \defeq
\{x\in \Mult(\coef\otimes\QG) \mid
x^*x\in \dom (\ID_\coef\otimes\Haar)^+\}.
\]

For any \(\psi\in\Sub(\Haar)\), there is a vector
\(\chi_\psi\in\Hilsr\) with
\begin{equation}
  \label{eq:def_cyclic_psi}
  \psi(x^*y) = \langle x\chi_\psi,y\chi_\psi\rangle
  \qquad\text{for all \(x,y\in\QG\),}
\end{equation}
and there is an operator~\(T_\psi\) on~\(\Hilsr\) with \(0\le
T_\psi\le 1\) that commutes with~\(\QG\) and satisfies
\begin{equation}
  \label{eq:cut-off_psi}
  \bigl\langle \Lambda(x),T_\psi\Lambda(y) \bigr\rangle
  = \psi(x^*y) \qquad\text{for all \(x,y\in\dom \Lambda\).}
\end{equation}
The net \((T_\psi)_{\psi\in\Sub(\Haar)}\) is increasing and
converges strongly towards \(\ID_\Hilsr\).  The
operator~\(T_\psi\) and the vector~\(\chi_\psi\) are related by
\begin{equation}
  \label{eq:cut-off_cyclic}
  T_\psi^{\nicefrac{1}{2}}\Lambda(x) = x\cdot\chi_\psi
  \qquad\text{for \(x\in\dom\Lambda\)}
\end{equation}
(see \cite{Kustermans-Vaes:LCQG}*{Notation 1.4}).

\begin{lemma}
  \label{lem:dom_slice_Lambda}
  We have \(x\in \dom(\ID_\coef\otimes\Lambda)\) if and only if
  the net \((x\cdot(1\otimes\chi_\psi))_{\psi\in\Sub(\Haar)}\) converges
  strongly in \(\Bound(\coef,\coef\otimes\Hilsr)\).
  Equivalently, \(x\cdot
  (b\otimes\chi_\psi)\in\coef\otimes\Hilsr\) converges in norm
  for each \(b\in\coef\).
\end{lemma}

This allows us to define \((\ID_\coef\otimes\Lambda)(x) \in
\Bound(\coef,\coef\otimes\Hilsr)\) by
\[
\bigl((\ID_\coef\otimes\Lambda)(x)\bigr)(b)
\defeq \lim_{\psi\in\Sub(\Haar)} x\cdot(b\otimes\chi_\psi)
\qquad\text{for all \(b\in\coef\).}
\]

\begin{proof}
  If
  \(\bigl(x\cdot(b\otimes\chi_\psi)\bigr)_{\psi\in\Sub(\Haar)}\)
  converges in norm for all \(b\in\coef\), then so does
  \[
  \bigl\langle x(b\otimes\chi_\psi),
  x(b\otimes\chi_\psi)\bigr\rangle_\coef =
  b^*(\ID\otimes\psi)(x^*x)b,
  \]
  and this is equivalent to the strict convergence of
  \((\ID\otimes\psi)(x^*x)\) by
  \cite{Kustermans-Vaes:Weight}*{Result 3.4}.  Hence
  \(x\in\dom(\ID_\coef\otimes\Lambda)\).

  Conversely, assume that
  \(x\in\dom(\ID_\coef\otimes\Lambda)\).  Let
  \(\ID_\coef\otimes\Lambda'\) be the map defined
  in~\cite{Kustermans-Vaes:Weight}.  Then
  \cite{Kustermans-Vaes:Weight}*{Result 3.21} yields
  \[
  (\ID\otimes T_\psi)\bigl(\ID\otimes\Lambda'(x)\bigr)\cdot b =
  x\cdot (b\otimes\chi_\psi)
  \]
  for all \(\psi\in\Sub(\Haar)\), \(b\in\coef\).  Since the net
  \((T_\psi)\) converges strongly to the identity map, so does
  \((\ID_\coef\otimes T_\psi)\).  Hence the net \(x\cdot
  (b\otimes\chi_\psi)\) converges towards
  \(\ID\otimes\Lambda'(x)\cdot b\).
\end{proof}

The proof also shows that our definition of
\(\ID\otimes\Lambda\) agrees with the one
in~\cite{Kustermans-Vaes:Weight}.  A direct proof that the net
\(\bigl(x\cdot (b\otimes\chi_\psi)\bigr)_{\psi\in\Sub(\Haar)}\)
in \(\coef\otimes\Hilsr\) converges for all
\(x\in\dom(\ID_\coef\otimes\Lambda)\) is similar to the last
part of the proof of Lemma~\ref{lem:coefficient_via_coaction}.

The above definition of \(\ID\otimes\Lambda\) shows that
\begin{equation}
  \label{eq:BRA_via_cut-off}
  \BRA{\xi}\eta
  \defeq (\ID\otimes\Lambda)
  \bigl(\coact_{\Mult(\Hilm)}(\xi)^* (\eta\otimes1_\QG)\bigr)
  = \lim_{\psi\in\Sub(\Haar)} \coact_{\Mult(\Hilm)}(\xi)^*
  (\eta\otimes\chi_\psi)
\end{equation}
for all \(\xi\in\Mult(\Hilm)_\si\).  It is clear from this that
the map \(\eta\mapsto \BRA{\xi}\eta\) is \(\coef\)\nb-linear.

\begin{lemma}
  \label{lem:BRA_not_multiplier}
  Let \(\xi\in\Mult(\Hilm)_\si\) and \(\eta\in\Hilm\).  Then
  \(\BRA{\xi}\eta\) belongs to \(\coef\otimes\Hilsr\) and not
  just to \(\Mult(\coef\otimes\Hilsr)\defeq
  \Bound(\coef,\coef\otimes\Hilsr)\).
\end{lemma}

\begin{proof}
  For any Hilbert \(\coef\)\nb-module, we have
  \(\Hilm\cdot\coef=\Hilm\).  Writing \(\eta=\eta'\cdot b\)
  with \(\eta'\in\Hilm\), \(b\in\coef\), we get \(\BRA{\xi}\eta
  = (\BRA{\xi}\eta')\cdot b\) by \(\coef\)\nb-linearity.
\end{proof}

\begin{example}
  \label{exa:si_Hilm_group}
  The definitions of square-integrable elements and the
  operator \(\BRA{\xi}\) above generalise the corresponding
  ones for group actions in \cites{Meyer:Equivariant,
    Meyer:Generalized_Fixed}.

  To see this, we mainly have to identify the map
  \(\ID_\coef\otimes\Lambda\) in the case \(\QG=\CONT_0(G)\)
  for a locally compact group~\(G\).  We describe
  \(\Sub(\Haar)\) as in Example~\ref{exa:integrable_group}.
  The map \(\ID_\coef\otimes\Lambda\) is defined on a bounded
  continuous function \(f\colon G\to\coef\) if and only if the
  net \((f\cdot w)_{w(g)\,\dd g\in\Sub(\Haar)}\) converges in
  \(L^2(G,\coef) \defeq \coef\otimes L^2(G)\).
  In~\cites{Meyer:Equivariant, Meyer:Generalized_Fixed}, a net
  of compactly supported continuous functions \(\chi_i\colon
  G\to[0,1]\) with \(\chi_i\to1\) uniformly on compact subsets
  of~\(G\) is used instead of the net of functions~\(w\) above.
  But it is not hard to see that both types of cut-off functions
  yield the same unbounded densely defined map
  \(\ID_\coef\otimes\Lambda\) from \(\CONT_b(G,\coef)\) to
  \(L^2(G,\coef)\).

  Furthermore, \(\coact_\Hilm(\xi)^*(\eta\otimes 1)\) for
  \(\xi,\eta\in\Hilm\) is the function \(g\mapsto \langle
  g\cdot \xi,\eta\rangle\).  Hence our definitions of the
  subspace~\(\Hilm_\si\) and the operator \(\BRA{\xi}\)
  specialise to the corresponding definitions in
  \cites{Meyer:Equivariant, Meyer:Generalized_Fixed} for
  \(\QG=\CONT_0(G)\).
\end{example}

\begin{example}
  \label{exa:si_Hilbert_space}
  If the coefficient algebra~\(\coef\) is trivial, then we are
  dealing with corepresentations of the quantum group
  \((\QG,\com)\) on Hilbert spaces.  In this case, our
  definition of a square-integrable vector agrees with the
  corresponding one in Section~\ref{sec:si_Hilbert_space} by
  Lemma~\ref{lem:coefficient_via_coaction}.  This provides some
  examples of square-integrable coactions on Hilbert modules.
  In particular, the left and right regular corepresentations
  on~\(\Hilsr\) are square-integrable by
  Lemma~\ref{lem:regular_si}.
\end{example}

\begin{example}
  \label{exa:compact_Hilm}
  Let \((\QG,\com)\) be a compact quantum group, so
  that~\(\QG\) is unital and the Haar weight~\(\Haar\) is
  bounded.  Then the maps \(\ID\otimes\Haar\) and
  \(\ID\otimes\Lambda\) are defined everywhere, so that every
  element of a Hilbert module is square-integrable.  Letting
  \(\Omega\defeq\Lambda(1_\QG)\in\Hilsr\), we have \(\Lambda(x)
  = x\cdot\Lambda(1_\QG) = x\cdot\Omega\) for all \(x\in\QG\).
  Hence
  \[
  \BRA{\xi}\eta
  = \coact_{\Mult(\Hilm)}(\xi)^*(\eta\otimes\Omega)
  \]
  for all \(\xi\in\Mult(\Hilm)\), \(\eta\in\Hilm\).
\end{example}

Our next goal is to show that the operator \(\BRA{\xi}\colon
\Hilm\to \coef\otimes\Hilsr\) for \(\xi\in\Mult(\Hilm)_\si\) is
adjointable and to describe its adjoint.  This involves the
left \(\QG^\dual\)\nb-module structure on \(\Mult(\Hilm)\)
defined by \(\omega\star\xi\defeq
(\ID\otimes\omega)\coact_{\Mult(\Hilm)}(\xi)\) for all
\(\omega\in \QG^\dual\), \(\xi\in\Mult(\Hilm)\).  Here we use
the slice map
\[
\ID\otimes\omega\colon \Bound(\coef\otimes\QG,\Hilm\otimes\QG)
\to \Bound(\coef,\Hilm)
\]
for a bounded weight~\(\omega\).

\begin{lemma}
  \label{lem:BRA_adjointable}
  If \(\xi\in\Mult(\Hilm)_\si\), then the operator
  \(\BRA{\xi}\colon \Hilm\to \coef\otimes\Hilsr\) is
  adjointable, and the adjoint operator \(\KET{\xi}\defeq
  \BRA{\xi}^*\) is determined by
  \[
  \KET{\xi}\bigl(b\otimes \Lambda(x)\bigr)
  = (\ID_\coef\otimes\Haar)
  \bigl(\coact_{\Mult(\Hilm)}(\xi)\cdot (b\otimes x)\bigr)
  = (x\Haar \star \xi)\cdot b
  \]
  for all \(b\in\coef\), \(x\in \dom \QGtodual\), that is,
  \(x\Haar\) is bounded and belongs to \(L^1(\QG)\).
\end{lemma}

\begin{proof}
  If \(\xi\in\Mult(\Hilm)_\si\), then the graph of the operator
  \(\BRA{\xi}\colon \Hilm\to \coef\otimes\Hilsr\) is closed
  because the unbounded operator \(\ID\otimes\Lambda\) is
  closed (see \cite{Kustermans-Vaes:Weight}*{Proposition
    3.23}).  The Closed Graph Theorem shows that \(\BRA{\xi}\)
  is bounded.

  Let \(b\in\coef\) and \(x\in\dom\QGtodual\).  Then
  \begin{multline*}
    \bigl\langle b\otimes \Lambda(x),\BRA{\xi}\eta\bigr\rangle
    = \bigl\langle b\otimes \Lambda(x),
    (\ID\otimes\Lambda)(\coact_{\Mult(\Hilm)}(\xi)^*)
    (\eta\otimes1)\bigr\rangle
    \\= (\ID\otimes\Haar)\bigl((b\otimes x)^*
    \coact_{\Mult(\Hilm)}(\xi)^*(\eta\otimes1)\bigr)
    = b^* (\ID\otimes\Haar)
    \bigl(\coact_{\Mult(\Hilm)}(\xi)(1\otimes x)\bigr)^*\eta
    \\= b^* \bigl((\ID\otimes x\Haar)
    \coact_{\Mult(\Hilm)}(\xi)\bigr)^*\eta
    = b^* (x\Haar\star\xi)^* \eta
    = \bigl\langle (x\Haar\star\xi)\cdot b,\eta\bigr\rangle.
  \end{multline*}
  Hence the adjoint of \(\BRA{\xi}\) is defined on
  \(b\otimes\Lambda(x)\) and maps it to \((x\Haar \star
  \xi)\cdot b\).  Since \(\coef\otimes \Lambda(\dom\QGtodual)\)
  is dense in \(\coef\otimes\Hilsr\) and \(\BRA{\xi}\) is bounded,
  this implies that \(\BRA{\xi}^*\) is defined on all of
  \(\coef\otimes\Hilsr\).
\end{proof}

For \(\psi\in\Sub(\Haar)\), we may use
\eqref{eq:def_cyclic_psi} and~\eqref{eq:cut-off_cyclic} to
rewrite
\begin{equation}
  \label{eq:psi_through_Lambda}
  \psi(y) = \langle \chi_\psi,y\cdot\chi_\psi\rangle
  = \bigl\langle \chi_\psi,
  T_\psi^{\nicefrac{1}{2}}\Lambda(y)\bigr\rangle
  = \bigl\langle T_\psi^{\nicefrac{1}{2}}\chi_\psi,
  \Lambda(y)\bigr\rangle.
\end{equation}
Hence \(\ID_\coef\otimes\psi\colon \coef\otimes\QG\to\coef\)
factors through the map \(\ID_\coef\otimes\Lambda\) as
\((\ID_\coef\otimes\psi)y = \bigl\langle 1\otimes
T_\psi^{\nicefrac{1}{2}}\chi_\psi,
(\ID_\coef\otimes\Lambda)(y)\bigr\rangle\); here we use
\(1\otimes T_\psi^{\nicefrac{1}{2}}\chi_\psi \in
\Mult(\coef\otimes\Hilsr) \defeq
\Bound(\coef,\coef\otimes\Hilsr)\), and the inner product
actually means that we apply the adjoint of this operator to
\((\ID_\coef\otimes\Lambda)(y)\in \coef\otimes\Hilsr\).  Using
this notation, we may write
\begin{equation}
  \label{eq:KET_without_splitting}
  \KET{\xi}(x)
  = \lim_{\psi\in\Sub(\Haar)}
  \bigl\langle 1\otimes T_\psi^{\nicefrac{1}{2}}\chi_\psi,
  \coact_{\Mult(\Hilm)}(\xi)(x)\bigr\rangle
  \qquad\text{for all \(x\in \coef\otimes\Hilsr\).}
\end{equation}
This formula works for all elements of \(\coef\otimes\Hilsr\),
not just for elementary tensors.

\begin{proposition}
  \label{pro:compose_BRA_KET}
  Let \(\xi,\eta\in\Mult(\Hilm)_\si\) and consider
  \(\xi\circ\eta^*\in\Bound(\Hilm)=\Mult\bigl(\Comp(\Hilm)\bigr)\).
  Then \(\xi\circ\eta^*\) is an integrable multiplier of
  \(\Comp(\Hilm)\) and \(\Av(\xi\circ\eta^*) =
  \KET{\xi}\BRA{\eta}\).  Conversely, if
  \(\xi\circ\xi^*\in\Mult\bigl(\Comp(\Hilm)\bigr)_\ible\) for
  some \(\xi\in\Mult(\Hilm)\), then~\(\xi\) is
  square-integrable.

  In particular, \(\ket{\xi}\bra{\eta}\in\Comp(\Hilm)_\ible\)
  for \(\xi,\eta\in\Hilm_\si\), and if \(\xi\in\Hilm\), then
  \(\xi\in\Hilm_\si\) if and only if
  \(\ket{\xi}\bra{\xi}\in\Comp(\Hilm)_\ible\).
\end{proposition}

\begin{proof}
  By polarisation, we may assume without loss of generality
  that \(\xi=\eta\).  We have to show that
  \(\xi\in\Mult(\Hilm)\) is square-integrable if and only if
  \(\xi\circ\xi^*\) is integrable and that
  \(\Av(\xi\circ\xi^*)=\KET{\xi}\BRA{\xi}\).  Recall
  that~\(\xi\) is square-integrable if and only if
  \(\coact_\Hilm(\xi)^*(\zeta\otimes1_\QG)\) belongs to the domain of
  \(\ID_\coef\otimes\Lambda\) for all \(\zeta\in\Hilm\).  By
  definition, this is equivalent to
  \((\zeta\otimes1_\QG)^*\coact_\Hilm(\xi)\coact_\Hilm(\xi)^*(\zeta\otimes1_\QG)\in
  \dom (\ID_\coef\otimes\Haar)^+\).  This is equivalent to
  strict convergence of
  \[
  (\ID_\coef\otimes\psi)\bigl((\zeta\otimes1_\QG)^*
  \coact_\Hilm(\xi)\coact_\Hilm(\xi)^*(\zeta\otimes1_\QG)\bigr)
  = \zeta^*\cdot
  (\ID_{\Bound(\Hilm)}\otimes\psi)\bigl(\coact_\Hilm(\xi\xi^*)\bigr)
  \cdot\zeta
  \]
  for \(\psi\in\Sub(\Haar)\), still for all \(\zeta\in\Hilm\).
  Using \cite{Kustermans-Vaes:Weight}*{Result 3.4}, this is
  equivalent to norm-convergence of \(b^*\zeta^*\cdot
  (\ID_{\Comp(\Hilm)}\otimes\psi)\bigl(\coact_\Hilm(\xi\xi^*)\bigr)
  \cdot\zeta b\) for all \(b\in\coef\), for all
  \(\zeta\in\Hilm\); using a variant of
  \cite{Kustermans-Vaes:Weight}*{Result 3.4}, whose proof is
  similar, this is equivalent to convergence of
  \((\ID_{\Comp(\Hilm)}\otimes\psi)\bigl(\coact_\Hilm(\xi\xi^*)\bigr)\)
  in the strict topology on \(\Bound(\Hilm) =
  \Mult\bigl(\Comp(\Hilm)\bigr)\).  By definition, this means
  that \(\xi\xi^*\) is integrable.  Thus~\(\xi\) is
  square-integrable if and only if \(\xi\xi^*\) is integrable.

  Our computation also shows that \(\bigl\langle \zeta,
  \KET{\xi}\BRA{\xi}\zeta\bigr\rangle_\coef = \langle \zeta,
  \Av(\xi\xi^*)\zeta\rangle_\coef\) for all \(\zeta\in\Hilm\)
  if~\(\xi\) is square-integrable.  This implies
  \(\KET{\xi}\BRA{\xi}=\Av(\xi\xi^*)\) by polarisation.  If
  \(\xi,\eta\in\Hilm\), then \(\xi\eta^*=\ket{\xi}\bra{\eta}\),
  so that we get the assertions in the second paragraph.
\end{proof}

\subsection{Conditions for square-integrability and examples}
\label{sec:si_conditions}

\begin{proposition}
  \label{pro:si_Comp_ible}
  A coaction on a Hilbert module~\(\Hilm\) is square-integrable
  if and only if the induced coaction on \(\Comp(\Hilm)\) is
  integrable.
\end{proposition}

\begin{proof}
  If~\(\Hilm\) is square-integrable, that is, \(\Hilm_\si\) is
  dense in~\(\Hilm\), then the elements \(\ket{\xi}\bra{\eta}\)
  for \(\xi,\eta\in\Hilm_\si\) span a dense subspace of
  \(\Comp(\Hilm)\).  Since these elements are integrable by
  Proposition~\ref{pro:compose_BRA_KET}, \(\Comp(\Hilm)_\ible\)
  is dense in \(\Comp(\Hilm)\), that is, the coaction on
  \(\Comp(\Hilm)\) is integrable (see
  Lemma~\ref{lem:integrable_coaction}).

  Conversely, assume that the coaction on \(\Comp(\Hilm)\) is
  integrable.  Then \(\Comp(\Hilm)\) contains a strictly
  positive integrable element~\(T\) by
  Lemma~\ref{lem:integrable_coaction}.  Since the positive
  integrable elements form a hereditary cone,
  \(T\ket{\xi}\bra{\xi}T = \ket{T\xi}\bra{T\xi}\) is integrable
  for all \(\xi\in\Hilm\), so that \(T\xi\) is
  square-integrable again by
  Proposition~\ref{pro:compose_BRA_KET}.  The strict positivity
  of~\(T\) means that the range of~\(T\) is dense in~\(\Hilm\).
  Since the range of~\(T\) consists of square-integrable
  elements, these are dense in~\(\Hilm\).
\end{proof}

\begin{example}
  \label{exa:si_on_Cstar}
  Let \(\Hilm=\coef\) viewed as a Hilbert \(\coef\)\nb-module.
  The given coaction on~\(\coef\) turns this into a
  \(\QG\)\nb-equivariant Hilbert \(\coef\)\nb-module.  Since
  the natural isomorphism \(\Comp(\coef)\cong\coef\) is
  \(\QG\)\nb-equivariant, Proposition~\ref{pro:si_Comp_ible}
  shows that~\(\coef\) is square-integrable as a Hilbert
  \(\coef\)\nb-module if and only if the coaction on the
  \(C^*\)\nb-algebra~\(\coef\) is integrable in the sense of
  Definition~\ref{def:integrable_coaction}.
\end{example}

Examples \ref{exa:si_Hilbert_space} and~\ref{exa:si_on_Cstar}
show that square-integrable coactions on Hilbert modules
contain both square-integrable Hilbert space corepresentations
and integrable coactions on \(C^*\)\nb-algebras as special
cases.

\begin{proposition}
  \label{pro:si_from_action}
  Let \(\coef\) and~\(\coefd\) be \(C^*\)\nb-algebras equipped
  with coactions of a locally compact quantum group
  \((\QG,\com)\), let~\(\Hilm\) be a \(\QG\)\nb-equivariant
  Hilbert \(\coef\)\nb-module, and let \(\pi\colon
  \coefd\to\Bound(\Hilm)\) be a \(\QG\)\nb-equivariant,
  non-degenerate \Star{}homomorphism.  If the coaction
  on~\(\coefd\) is integrable, then~\(\Hilm\) is
  square-integrable.
\end{proposition}

\begin{proof}
  The coaction on \(\Comp(\Hilm)\) is integrable by
  Proposition~\ref{pro:integrable_functorial}.  This is
  equivalent to square-integrability of~\(\Hilm\) by
  Proposition~\ref{pro:si_Comp_ible}.
\end{proof}

As in Proposition~\ref{pro:integrable_functorial}, we may
replace~\(\pi\) by a positive, strictly continuous,
non-degenerate linear map \(\Mult(\coefd)\to\Bound(\Hilm)\).

\begin{proposition}
  \label{pro:si_stable}
  Let~\(\coef\) be a \(C^*\)\nb-algebra equipped with a
  coaction of a locally compact quantum group.  Equip
  \(\coef\otimes\Hilsr\) with the diagonal coaction defined
  in~\eqref{eq:standard_regular_coact}.  Then
  \(\coef\otimes\Hilsr\) is square-integrable.  More generally,
  if~\(\Hilm\) is a \(\QG\)\nb-equivariant Hilbert
  \(\coef\)\nb-module, then \(\Hilm\otimes\Hilsr\) with the
  diagonal coaction is square-integrable.
\end{proposition}

\begin{proof}
  Tensor product coactions on Hilbert modules are defined
  in~\cite{Baaj-Skandalis:Hopf_KK}.  This specialises
  to~\eqref{eq:standard_regular_coact} for
  \(\coef\otimes\Hilsr\), so that the latter is a
  \(\QG\)\nb-equivariant Hilbert \(\coef\)\nb-module.
  Proposition~\ref{pro:stable_integrable} asserts that
  \(\Comp(\Hilm\otimes\Hilsr) \cong
  \Comp\bigl(\Comp(\Hilm)\otimes\Hilsr\bigr)\) is integrable.
  Hence \(\Hilm\otimes\Hilsr\) is square-integrable by
  Proposition~\ref{pro:si_Comp_ible}.
\end{proof}

\begin{proposition}
  \label{pro:si_sum}
  Let \((\Hilm_j)_{j\in I}\) be a set of \(\QG\)\nb-equivariant
  Hilbert \(\coef\)\nb-modules.  Let \(\Hilm'\defeq
  \bigoplus_{j\in I} \Hilm_j\) be the Hilbert module completion
  of the algebraic direct sum of the
  \(\coef\)\nb-modules~\(\Hilm_j\) for the obvious inner
  product \(\bigl\langle (\xi_j),(\eta_j) \bigr\rangle \defeq
  \sum_{j\in I} \langle\xi_j,\eta_j\rangle\)

  The induced coaction of~\(\QG\) on~\(\Hilm'\) is
  square-integrable if and only if~\(\Hilm_j\) is
  square-integrable for each \(j\in I\).
\end{proposition}

\begin{proof}
  The \(C^*\)\nb-algebra \(\Comp(\Hilm_j)\) for some \(j\in I\)
  is a corner in \(\Comp(\Hilm')\).  Using
  Proposition~\ref{pro:si_Comp_ible} and
  Corollary~\ref{cor:corner_integrable}, we see
  that~\(\Hilm_j\) must be square-integrable if~\(\Hilm'\) is.
  Conversely, it is straightforward to see that a vector
  \((\xi_j)_{j\in I}\) in \(\Hilm'\defeq \bigoplus \Hilm_j\)
  with \(\xi_j=0\) for all but finitely many  \(j\in I\) is
  square-integrable if and only if \(\xi_j\in\Hilm_j\) is
  square-integrable for all \(j\in I\).  Hence the
  square-integrable elements are dense in~\(\Hilm'\) if this
  holds for~\(\Hilm_j\) for all \(j\in I\).
\end{proof}

For a finite direct sum \(\Hilm_1\oplus\dotsb\oplus\Hilm_n\),
it is easy to see that an element \((\xi_1,\dotsc,\xi_n)\) is
square-integrable if and only if the components
\(\xi_j\in\Hilm_j\) are square-integrable for \(j=1,\dotsc,n\).

It is unclear whether square-integrability is inherited by
Hilbert submodules that are not complemented.  This corresponds
to the question whether integrability passes on to hereditary
\(C^*\)\nb-subalgebras.

\subsection{Further properties of square-integrable elements}
\label{sec:si_further}

\begin{proposition}
  \label{pro:BRA_equivariant}
  The operators \(\BRA{\xi}\colon \Hilm\to \coef\otimes\Hilsr\)
  and \(\KET{\xi}\colon \coef\otimes\Hilsr\to\Hilm\) for
  square-integrable~\(\xi\) are \(\QG\)\nb-equivariant.
\end{proposition}

Here we use the coaction on \(\coef\otimes\Hilsr\) defined
in~\eqref{eq:standard_regular_coact}.

\begin{proof}
  A direct proof of this assertion is possible but somewhat
  unpleasant because the coaction on \(\coef\otimes\Hilsr\) is
  complicated.  We avoid this direct proof by a trick.  Recall
  that \(\KET{\xi}\BRA{\eta} = \Av(\xi^*\eta)\) for all
  \(\xi,\eta\in\Mult(\Hilm)\) by
  Proposition~\ref{pro:compose_BRA_KET}.  This is
  \(\QG\)\nb-equivariant by Lemma~\ref{lem:av_invariant}.
  Replacing \(\Hilm\) by \(\Hilm_1\oplus\Hilm_2\), we see that
  the same holds for \(\KET{\xi}\BRA{\eta}\colon
  \Hilm_2\to\Hilm_1\) for \(\xi\in\Mult(\Hilm_1)_\si\) and
  \(\eta\in\Mult(\Hilm_2)_\si\).  Now we consider
  \(\Hilm_2\defeq \coef\otimes\Hilsr\) with the diagonal
  coaction and use a multiplier of the form
  \(1_\coef\otimes\eta\) with \(\eta\in\Hilsr_\si\).  It is
  easy to check that \(\BRA{1_\coef\otimes\eta} =
  1_\coef\otimes\BRA{\eta}\) in this case, and the equivariance
  of \(\BRA{\eta}\colon \Hilsr\to\Hilsr\) is already checked in
  Lemma~\ref{lem:KET_equivariant}.

  Summing up, if \(\xi\in\Mult(\Hilm)_\si\), then we know that
  the operators
  \[
  \KET{\xi}\BRA{1_\coef\otimes\eta}\colon
  \coef\otimes\Hilsr\to\Hilm
  \qquad\text{and}\qquad
  \BRA{1_\coef\otimes\eta}\colon
  \coef\otimes\Hilsr\to\coef\otimes\Hilsr
  \]
  are \(\QG\)\nb-equivariant for all \(\eta\in\Hilsr_\si\).
  Thus
  \begin{multline*}
    \coact_\Hilm\bigl(\KET{\xi}\BRA{1_\coef\otimes\eta}h\bigr)
    = \bigl(\KET{\xi}\BRA{1_\coef\otimes\eta}\otimes 1_\QG\bigr)
    \bigl((\coact_\coef\otimes\ID_\Hilsr)(h)\bigr)
    \\= \bigl(\KET{\xi}\otimes 1_\QG\bigr)
    \bigl(\coact_\Hilm(\BRA{1_\coef\otimes\eta}h)\bigr)
  \end{multline*}
  for all \(h\in\coef\otimes\Hilsr\).

  We claim that the sum of the ranges of the operators
  \(\BRA{1_\coef\otimes\eta}=1_\coef\otimes\BRA{\eta}\) for
  \(\eta\in\Hilsr_\si\) is dense in \(\coef\otimes\Hilsr\).  Hence
  \[
  \coact_\Hilm(\KET{\xi}k)
  = (\KET{\xi}\otimes 1_\QG) \coact_{\coef\otimes\Hilsr}(k)
  \]
  holds for a dense set of \(k\in\coef\otimes\Hilsr\) and hence
  for all \(k\in\coef\otimes\Hilsr\) by continuity.  This shows
  that \(\KET{\xi}\) is \(\QG\)\nb-equivariant.  Hence so is
  its adjoint \(\BRA{\xi}\).

  It remains to check the claim above.  It suffices to prove
  that the ranges of the operators \(\BRA{\eta}\) for
  \(\eta\in\Hilsr_\si\) span a dense subspace in~\(\Hilsr\).
  This follows easily from the non-degeneracy of the right
  regular representation of~\(\hat{\QG}\) on~\(\Hilsr\) and the
  computations in the proof of Lemma~\ref{lem:regular_si}.
\end{proof}

\begin{lemma}
  \label{lem:equivariant_on_si}
  If \(T\in\Bound(\Hilm)\) is \(\QG\)\nb-equivariant and
  \(\xi\in\Mult(\Hilm)_\si\), then
  \[
  T\xi\in\Mult(\Hilm)_\si\qquad\text{with}\quad
  \BRA{T\xi} = \BRA{\xi}\circ T^*,\quad
  \KET{T\xi} = T\circ\KET{\xi}.
  \]
\end{lemma}

\begin{proof}
  Since~\(T\) is equivariant,
  \(\coact_\Hilm(T\xi)=(T\otimes1_\QG)\coact_\Hilm(\xi)\), so that we may
  rewrite \(\coact_\Hilm(T\xi)^*(\eta\otimes1_\QG) =
  \coact_\Hilm(\xi)^*(T^*\eta\otimes1_\QG)\).  This implies the
  assertions.
\end{proof}

\begin{lemma}
  \label{si_complete}
  The subspaces \(\Mult(\Hilm)_\si\subseteq\Mult(\Hilm)\) and
  \(\Hilm_\si\subseteq\Hilm\) are both complete with respect to
  the norm
  \[
  \norm{\xi}_\si \defeq \norm{\xi} +
  \bigl\lVert\KET{\xi}\bigr\rVert.
  \]
\end{lemma}

\begin{proof}
  This is the graph norm for the unbounded linear map
  \(\xi\mapsto \KET{\xi}\).  Hence the assertion boils down to
  the claim that this operator with domain \(\Mult(\Hilm)_\si\)
  or \(\Hilm_\si=\Mult(\Hilm)_\si\cap\Hilm\) is closed.  This
  is true because \(\coact_{\Mult(\Hilm)}\) is bounded and
  \(\ID_\coef\otimes\Lambda\) is closed by
  \cite{Kustermans:KMS}*{Result 2.3}.
\end{proof}

\begin{lemma}
  \label{lem:si_multiply}
  Let \(\xi\in\Mult(\Hilm)_\si\) and \(b\in\Mult(\coef)\). Then
  \[
  \xi\cdot b\in\Mult(\Hilm)_\si\qquad\text{with}\quad
  \BRA{\xi\cdot b} = \coact_\coef(b)^* \circ \BRA{\xi},\quad
  \KET{\xi\cdot b} = \KET{\xi}\circ \coact_\coef(b),
  \]
  where we represent
  \(\coact_\coef(b)\in\Mult(\coef\otimes\QG)\) on
  \(\coef\otimes\Hilsr\) using the standard representation
  of~\(\QG\) on~\(\Hilsr\).
\end{lemma}

\begin{proof}
  The subspace \(\dom(\ID_\coef\otimes\Lambda)\) is a left
  ideal in \(\Mult(\coef\otimes\QG)\) because the subset
  \(\dom(\ID_\coef\otimes\Haar)^+\) is a hereditary cone (see
  also~\cite{Kustermans-Vaes:Weight}).  Furthermore,
  \(\ID_\coef\otimes\Lambda\) is a left module homomorphism.
  Using \(\coact_{\Hilm}(\xi\cdot b)^* =
  \coact_{\coef}(b)^* \coact_{\Hilm}(\xi)^*\), we
  get \(\xi\cdot b\in\Mult(\Hilm)_\si\) and
  \begin{multline*}
    \BRA{\xi\cdot b}\eta
    = (\ID_\coef\otimes\Lambda)
    \bigl(\coact_{\Hilm}(\xi b)^*\cdot (\eta\otimes1_\QG)\bigr)
    = (\ID_\coef\otimes\Lambda)
    \bigl(\coact_{\coef}(b)^*\cdot\coact_{\Hilm}(\xi)^*\cdot
    (\eta\otimes1_\QG)\bigr)
    \\= \coact_{\coef}(b)^* (\ID_\coef\otimes\Lambda)
    \bigl(\coact_{\Hilm}(\xi)^*\cdot (\eta\otimes1_\QG)\bigr)
    = \coact_{\coef}(b)^* \BRA{\xi}\eta
  \end{multline*}
  for all \(\eta\in\Hilm\).  Hence \(\BRA{\xi\cdot
    b}=\coact_{\coef}(b)^*\BRA{\xi}\).  The formula for \(\KET{\xi\cdot
    b}\) follows by taking adjoints.
\end{proof}

Lemma~\ref{lem:si_multiply} can be generalised as follows.  Let
\(\Hilm_1\) and~\(\Hilm_2\) be two \(\QG\)\nb-equivariant
Hilbert \(\coef\)\nb-modules, let
\(T\in\Bound(\Hilm_1,\Hilm_2)\) and
\(\xi\in\Mult(\Hilm_1)\defeq \Bound(\coef,\Hilm_1)\).  All
relevant Hilbert \(\coef\)\nb-modules embed in \(\Hilm\defeq
\Hilm_1\oplus\Hilm_2\oplus\coef\).  We may view~\(T\)
and~\(\xi\) as multipliers of \(\Comp(\Hilm)\).  If~\(T\) is
square-integrable as such a multiplier, then
Lemma~\ref{lem:si_multiply} yields that \(T\xi\) is
square-integrable as well and satisfies
\[
\BRA{T\xi} = \coact_{\Hilm_1}(\xi)^* \BRA{T},\qquad
\KET{T\xi} = \KET{T}\coact_{\Hilm_1}(\xi).
\]
These formulas originally involve two operators
\(\Comp(\Hilm)\otimes\Hilsr\to\Comp(\Hilm)\).  But it continues
to hold if we reinterpret \(\KET{T}\) as an operator
\(\Hilm_1\otimes\Hilsr\to\Hilm_2\), \(\coact_{\Hilm_1}(\xi)\) as an
operator \(\coef\otimes\Hilsr\to\Hilm_1\otimes\Hilsr\), and
\(\KET{T\xi}\) as an operator \(\coef\otimes\Hilsr\to\Hilm_2\).

The case \(\Hilm_1=\Hilm_2\) is particularly interesting
because we may get square-integrable operators
\(\Hilm_1\to\Hilm_1\) from an equivariant \Star{}representation
\(\coefd\to\Bound(\Hilm_1)\) and square-integrable multipliers
of~\(\coefd\).  As a consequence,
\(\Mult(\coefd)_\si\cdot\Mult(\Hilm)\subseteq\Mult(\Hilm)_\si\).
This yields an alternative proof of
Proposition~\ref{pro:si_stable} that only uses
square-integrable elements, not integrable elements.

Let \(\xi\in\Mult(\Hilm)_\si\) and
\(\omega\in\QG^\dual\subseteq \Mult(\hat{\QG})\).  We ask when
\(\omega\star\xi\) is square-integrable as well and how to
describe \(\KET{\omega\star\xi}\).  This requires the modular
automorphism group \(\hat{\modauto}\colon
\R\times\hat{\QG}\to\hat{\QG}\) and the right regular
anti-representation defined in~\eqref{eq:right_regular_anti}.

\begin{lemma}
  \label{lem:si_act}
  Let \(\xi\in\Mult(\Hilm)_\si\) and
  \(\omega\in\QG^\dual\subseteq \Mult(\hat{\QG})\).  If
  \(\omega\in \dom(\hat{\modauto}_{\nicefrac{\ima}{2}})\), then
  \[
  \omega\star\xi\in\Mult(\Hilm)_\si
  \qquad\text{with}\quad
  \KET{\omega\star\xi} = \KET{\xi}\circ \bigl(\ID_\coef\otimes
  \rho\circ\hat{\modauto}_{\nicefrac{\ima}{2}}(\omega)\bigr).
  \]
\end{lemma}

\begin{proof}
  Recall that \(\hat{\Lambda}(x\Haar) = \Lambda(x)\) for
  \(x\in\dom\QGtodual\).  Hence we may rewrite the formula for
  \(\KET{\xi}\) in Lemma~\ref{lem:BRA_adjointable} as
  \(\KET{\xi}\bigl(b\otimes\hat{\Lambda}(x)\bigr) =
  (x\star\xi)\cdot b\) for all \(b\in\coef\), \(x\in
  \QG^\dual\cap\dom(\hat{\Lambda})\).  For these \(b\)
  and~\(x\), the formula makes sense regardless whether~\(\xi\)
  is square-integrable or not.  Hence we may compute
  \(\KET{\omega\star\xi}\bigl(b\otimes\hat{\Lambda}(x)\bigr)\).

  Since the product~\(\star\) is associative, we get
  \[
  \KET{\omega\star\xi}\bigl(b\otimes\hat{\Lambda}(x)\bigr)
  = (x\star\omega\star\xi)\cdot b
  = \KET{\xi}\bigl(b\otimes \hat{\Lambda}(x\omega)\bigr).
  \]
  Now we use~\eqref{eq:right_regular_rep} to rewrite the right
  hand side as
  \[
  \KET{\xi}\bigl(b\otimes
  \rho\circ\hat{\modauto}_{\nicefrac{\ima}{2}}(\omega)
  \hat{\Lambda}(x)\bigr)
  = \KET{\xi}\circ \bigr(\ID_\coef\otimes
  \rho\circ\hat{\modauto}_{\nicefrac{\ima}{2}}(\omega)\bigr)
  \bigl(b\otimes\hat{\Lambda}(x)\bigr).
  \]
  This yields the asserted formula for \(\KET{\omega\star\xi}\)
  and shows that \(\KET{\omega\star\xi}\) extends to an
  adjointable operator \(\coef\otimes\Hilsr\to\Hilm\).  The
  adjoint operator must be \(\BRA{\omega\star\xi}\), forcing
  \(\omega\star\xi\) to be square-integrable.
\end{proof}

\section{The Stabilisation Theorem}
\label{sec:stabilisation}

Now we come to our main result, an equivariant generalisation
of Kasparov's Stabilisation Theorem for actions of locally
compact quantum groups.

\begin{theorem}
  \label{the:stabilisation}
  Let \((\QG,\com)\) be a locally compact quantum group,
  let~\(\coef\) be a \(C^*\)\nb-algebra with a coaction of
  \((\QG,\com)\), and let~\(\Hilm\) be a countably generated
  \(\QG\)\nb-equivariant Hilbert \(\coef\)\nb-module.  Then the
  following are equivalent:
  \begin{enumerate}[label=\textup{(\arabic{*})}]
  \item there is a \(\QG\)\nb-equivariant unitary
    \(\Hilm\oplus \coef\otimes\Hilsr^\infty \cong
    \coef\otimes \Hilsr^\infty\); here
    \(\coef\otimes\Hilsr^\infty\) denotes the direct sum of
    countably many copies of the Hilbert \(\coef\)\nb-module
    \(\coef\otimes\Hilsr\), and the latter is equipped with the
    coaction of~\(\QG\) described
    in~\textup{\eqref{eq:standard_regular_coact}}.

  \item There is an adjointable \(\QG\)\nb-equivariant isometry
    \(\Hilm\to \coef\otimes\Hilsr^\infty\), that is,
    \(\Hilm\) is \(\QG\)\nb-equivariantly isomorphic to a
    direct summand of \(\coef\otimes\Hilsr^\infty\).

  \item The subspace of square-integrable elements is dense
    in~\(\Hilm\), that is, the coaction on~\(\Hilm\) is
    square-integrable.

  \item The subspace of integrable operators is dense in
    \(\Comp(\Hilm)\), that is, the induced coaction on
    \(\Comp(\Hilm)\) is integrable.
  \end{enumerate}
\end{theorem}

For a compact quantum group, \(\Hilm\oplus
\coef\otimes\Hilsr^\infty \cong \coef\otimes\Hilsr^\infty\)
holds for any countably generated equivariant Hilbert
\(\coef\)\nb-module~\(\Hilm\) because (3) and~(4) are
automatically true.

The following technical result is used in the proof of
Theorem~\ref{the:stabilisation}.

\begin{lemma}
  \label{lem:non-degenerate_module}
  For any \(\QG\)\nb-equivariant Hilbert
  \(\coef\)\nb-module~\(\Hilm\), we have
  \[
  \clspan \{(\omega\star\xi)\cdot b \mid
  \omega\in L^1(\QG),\ \xi\in\Hilm,\ b\in\coef\}
  = \Hilm.
  \]
\end{lemma}

\begin{proof}
  By definition, \((a\omega\star\xi)\cdot b =
  (\ID\otimes\omega) \bigl(\coact_\Hilm(\xi)\cdot(b\otimes
  a)\bigr)\) for \(a\in\QG\), \(\omega\in L^1(\QG)\),
  \(\xi\in\Hilm\), \(b\in\coef\).  The definition of a Hilbert
  module coaction contains the requirement that the elements
  \(\coact_\Hilm(\xi)\cdot(b\otimes a)\) span a dense subspace
  of \(\Hilm\otimes\QG\).  Applying \(\ID\otimes\omega\), we
  get a dense subspace of~\(\Hilm\).
\end{proof}

\begin{proof}[Proof of Theorem~\textup{\ref{the:stabilisation}}]
  The implication (1)\(\Longrightarrow\)(2) is trivial, and the
  equivalence of (3) and~(4) is already contained in
  Proposition~\ref{pro:si_Comp_ible}.

  We check (2)\(\Longrightarrow\)(3).
  Proposition~\ref{pro:si_stable} shows that
  \(\coef\otimes\Hilsr^\infty \cong \coef^\infty\otimes\Hilsr\)
  is square-integrable.  Direct summands of square-integrable
  equivariant Hilbert module inherit square-integrability by
  Proposition~\ref{pro:si_sum}.  Hence any direct summand of
  \(\coef\otimes\Hilsr^\infty\) is square-integrable.  Thus~(2)
  implies~(3).

  Finally, we check that~(3) implies~(1).  If \(\xi\in\Hilm\)
  is square-integrable, then \(\KET{\xi}\) is an equivariant
  adjointable map \(\coef\otimes\Hilsr\to\Hilm\).  The formula
  for \(\KET{\xi}\) in Lemma~\ref{lem:BRA_adjointable} shows
  that its range contains \((x\cdot\Haar\star\xi)\cdot b\) for
  all \(x\in\dom\QGtodual\), \(b\in\coef\).  Since the set of
  \(x\cdot\Haar\) with \(x\in\dom\QGtodual\) is dense in
  \(L^1(\QG)\) and the set of square-integrable vectors is
  dense in~\(\Hilm\) by assumption,
  Lemma~\ref{lem:non-degenerate_module} shows that the ranges
  of the operators \(\KET{\xi}\) span a dense subspace
  of~\(\Hilm\).  This implies~(1) using a well-known argument
  by Mingo and Phillips (see~\cite{Mingo-Phillips:Triviality}),
  which is also used in~\cite{Meyer:Equivariant} to prove the
  Equivariant Stabilisation Theorem for Hilbert modules with an
  action of a locally compact group.  The idea of this argument
  is to construct an equivariant adjointable map
  \(\coef\otimes\Hilsr^\infty \to \Hilm\oplus
  \coef\otimes\Hilsr^\infty\) which is injective and has dense
  range.  Then a polar decomposition provides the desired
  equivariant unitary.  We refer to
  \cites{Mingo-Phillips:Triviality, Meyer:Equivariant} for
  further details.
\end{proof}

\begin{proposition}
  \label{pro:integrable_via_cutoff}
  A countably generated \(\QG\)\nb-equivariant Hilbert
  module~\(\Hilm\) is square-integrable if and only if there is
  a non-degenerate, equivariant, completely positive
  contractive linear map \(\QG\to\Bound(\Hilm)\).

  A coaction of~\(\QG\) on a \(\sigma\)\nb-unital
  \(C^*\)\nb-algebra~\(\coef\) is integrable if and only if
  there is a non-degenerate, equivariant, completely positive
  contractive linear map \(\QG\to\Mult(\coef)\).
\end{proposition}

\begin{proof}
  The two assertions are equivalent by
  Proposition~\ref{pro:si_Comp_ible}.
  Proposition~\ref{pro:integrable_functorial} shows that the
  existence of a map \(\QG\to\Mult(\coef)\) as in the theorem
  is sufficient for integrability.  Conversely, assume
  that~\(\Hilm\) is a  countably generated, square-integrable
  \(\QG\)\nb-equivariant Hilbert \(\coef\)\nb-module.   The
  Equivariant Stabilisation Theorem~\ref{the:stabilisation}
  provides an equivariant adjointable isometry \(T\colon
  \Hilm\to \coef\otimes\Hilsr^\infty\).  There are equivariant
  \Star{}homomorphisms \(\pi\colon\QG\to\Bound(\Hilsr)\to
  \Bound(\coef\otimes\Hilsr^\infty)\).  Then the map
  \[
  \QG\to\Bound(\Hilm),
  \qquad x\mapsto T^*\pi(x)T
  \]
  has the required properties.
\end{proof}

\section{Conclusion}
\label{sec:conclusion}

We have extended the theory of square-integrable group actions
on Hilbert modules to coactions of locally compact quantum
groups.  This includes basic results on irreducible
square-integrable Hilbert space corepresentations.  Our main
result, the Equivariant Stabilisation Theorem shows that
square-integrable coactions on Hilbert modules are closely
related to the regular corepresentation and hence to stable
coactions: a coaction on a Hilbert \(\coef\)\nb-module is
square-integrable if and only if it is contained in the
diagonal coaction on \(\coef\otimes\Hilsr^\infty\),
where~\(\Hilsr\) carries the left regular representation.

It is remarkable that we do not need the quantum group to be
regular or the coaction to be continuous (continuity of
coactions is only problematic for non-regular quantum groups,
see~\cite{Baaj-Skandalis-Vaes:Non-semi-regular}).

In the theory of square-integrable group actions, the next step
brings in crossed products.  In the quantum group setting, the
issue is whether or not the operators
\(\BRA{\xi}\circ\KET{\eta}\) on \(\coef\otimes\Hilsr\) belong
to the canonical representation of the crossed product
\(C^*\)\nb-algebra for sufficiently many square-integrable
\(\xi,\eta\).  If this is the case, we may define a generalised
fixed point algebra that is Morita equivalent to an ideal in
the crossed product.  As in~\cite{Meyer:Generalized_Fixed},
this provides an equivalence between Hilbert modules over the
crossed products and \(\QG\)\nb-equivariant Hilbert modules
over~\(\coef\) with a suitable dense subspace.  This line of
thought is pursued in the first author's doctoral thesis and
will be published elsewhere.


\begin{bibdiv}
\begin{biblist}
\bib{Baaj-Skandalis:Hopf_KK}{article}{
  author={Baaj, Saad},
  author={Skandalis, Georges},
  title={\(C^*\)\nobreakdash -alg\`ebres de Hopf et th\'eorie de Kasparov \'equivariante},
  language={French, with English summary},
  journal={\(K\)-Theory},
  volume={2},
  date={1989},
  number={6},
  pages={683--721},
  issn={0920-3036},
  review={\MRref {1010978}{90j:46061}},
}

\bib{Baaj-Skandalis:Unitaires}{article}{
  author={Baaj, Saad},
  author={Skandalis, Georges},
  title={Unitaires multiplicatifs et dualit\'e pour les produits crois\'es de $C^*$-alg\`ebres},
  language={French, with English summary},
  journal={Ann. Sci. \'Ecole Norm. Sup. (4)},
  volume={26},
  date={1993},
  number={4},
  pages={425--488},
  issn={0012-9593},
  review={\MRref {1235438}{94e:46127}},
}

\bib{Baaj-Skandalis-Vaes:Non-semi-regular}{article}{
  author={Baaj, Saad},
  author={Skandalis, Georges},
  author={Vaes, Stefaan},
  title={Non-semi-regular quantum groups coming from number theory},
  journal={Comm. Math. Phys.},
  volume={235},
  date={2003},
  number={1},
  pages={139--167},
  issn={0010-3616},
  review={\MRref {1969723}{2004g:46083}},
}

\bib{Combes:Poids}{article}{
  author={Combes, Fran\c {c}ois},
  title={Poids sur une $C^*$-alg\`ebre},
  language={French},
  journal={J. Math. Pures Appl. (9)},
  volume={47},
  date={1968},
  pages={57--100},
  issn={0021-7824},
  review={\MRref {0236721}{38 \#5016}},
}

\bib{Combes:Systemes_Hilbertiens}{article}{
  author={Combes, Fran\c {c}ois},
  title={Syst\`emes hilbertiens \`a gauche et repr\'esentation de Gelfand-Segal},
  language={French, with English summary},
  conference={ title={Operator algebras and group representations, Vol. I (Neptun, 1980)}, },
  book={ series={Monogr. Stud. Math.}, volume={17}, publisher={Pitman}, place={Boston, MA}, },
  date={1984},
  pages={71--107},
  review={\MRref {731765}{85f:46103}},
}

\bib{Exel:Unconditional}{article}{
  author={Exel, Ruy},
  title={Unconditional integrability for dual actions},
  journal={Bol. Soc. Brasil. Mat. (N.S.)},
  volume={30},
  number={1},
  date={1999},
  pages={99--124},
  issn={0100-3569},
  review={\MRref {1686980}{2000f:46071}},
}

\bib{Kustermans:KMS}{article}{
  author={Kustermans, Johan},
  title={KMS weights on \(C^*\)\nobreakdash -algebras},
  status={eprint},
  note={\arxiv {math/9704008}},
  date={1997},
}

\bib{Kustermans:LCQG_universal}{article}{
  author={Kustermans, Johan},
  title={Locally compact quantum groups in the universal setting},
  journal={Internat. J. Math.},
  volume={12},
  date={2001},
  number={3},
  pages={289--338},
  issn={0129-167X},
  review={\MRref {1841517}{2002m:46108}},
}

\bib{Kustermans-Vaes:Weight}{article}{
  author={Kustermans, Johan},
  author={Vaes, Stefaan},
  title={Weight theory for \(C^*\)\nobreakdash -algebraic quantum groups},
  status={eprint},
  note={\arxiv {math/9901063}},
  date={1999},
}

\bib{Kustermans-Vaes:LCQG}{article}{
  author={Kustermans, Johan},
  author={Vaes, Stefaan},
  title={Locally compact quantum groups},
  language={English, with English and French summaries},
  journal={Ann. Sci. \'Ecole Norm. Sup. (4)},
  volume={33},
  date={2000},
  number={6},
  pages={837--934},
  issn={0012-9593},
  review={\MRref {1832993}{2002f:46108}},
}

\bib{Kustermans-Vaes:LCQGvN}{article}{
  author={Kustermans, Johan},
  author={Vaes, Stefaan},
  title={Locally compact quantum groups in the von Neumann algebraic setting},
  journal={Math. Scand.},
  volume={92},
  date={2003},
  number={1},
  pages={68--92},
  issn={0025-5521},
  review={\MRref {1951446}{2003k:46081}},
}

\bib{Meyer:Equivariant}{article}{
  author={Meyer, Ralf},
  title={Equivariant Kasparov theory and generalized homomorphisms},
  journal={\(K\)\nobreakdash -Theory},
  volume={21},
  date={2000},
  number={3},
  pages={201--228},
  issn={0920-3036},
  review={\MRref {1803228}{2001m:19013}},
}

\bib{Meyer:Generalized_Fixed}{article}{
  author={Meyer, Ralf},
  title={Generalized fixed point algebras and square-integrable groups actions},
  journal={J. Funct. Anal.},
  volume={186},
  date={2001},
  number={1},
  pages={167--195},
  issn={0022-1236},
  review={\MRref {1863296}{2002j:46086}},
}

\bib{Mingo-Phillips:Triviality}{article}{
  author={Mingo, James A.},
  author={Phillips, William J.},
  title={Equivariant triviality theorems for Hilbert \(C^*\)\nobreakdash -modules},
  journal={Proc. Amer. Math. Soc.},
  volume={91},
  year={1984},
  number={2},
  pages={225--230},
  issn={0002-9939},
  review={\MRref {740176}{85f:46111}},
}

\bib{Rieffel:Integrable_proper}{article}{
  author={Rieffel, Marc A.},
  title={Integrable and proper actions on $C^*$-algebras, and square-integrable representations of groups},
  journal={Expo. Math.},
  volume={22},
  date={2004},
  number={1},
  pages={1--53},
  issn={0723-0869},
  review={\MRref {2166968}{2006g:46108}},
}

\bib{Stratila-Zsido:von_Neumann}{book}{
  author={Str\u {a}til\u {a}, \c {S}erban},
  author={Zsid\'o, L{\'a}szl\'o},
  title={Lectures on von Neumann algebras},
  note={Revision of the 1975 original; translated from the Romanian by Silviu Teleman},
  publisher={Editura Academiei},
  place={Bucharest},
  date={1979},
  pages={478},
  isbn={0-85626-109-2},
  review={\MRref {526399}{81j:46089}},
}

\bib{Takesaki:Theory_1}{book}{
  author={Takesaki, Masamichi},
  title={Theory of operator algebras. I},
  series={Encyclopaedia of Mathematical Sciences},
  volume={124},
  note={Reprint of the first (1979) edition; Operator Algebras and Non-commutative Geometry, 5},
  publisher={Springer-Verlag},
  place={Berlin},
  date={2002},
  pages={xx+415},
  isbn={3-540-42248-X},
  review={\MRref {1873025}{2002m:46083}},
}

\bib{Vaes:Induction_Imprimitivity}{article}{
  author={Vaes, Stefaan},
  title={A new approach to induction and imprimitivity results},
  journal={J. Funct. Anal.},
  volume={229},
  year={2005},
  number={2},
  pages={317--374},
  issn={0022-1236},
  review={\MRref {2182592}{2007f:46065}},
}

\bib{Woronowicz:Compact_pseudogroups}{article}{
  author={Woronowicz, Stanis{\l }aw Lech},
  title={Compact matrix pseudogroups},
  journal={Comm. Math. Phys.},
  volume={111},
  date={1987},
  number={4},
  pages={613--665},
  issn={0010-3616},
  review={\MRref {901157}{88m:46079}},
}

\bib{Woronowicz:CQG}{article}{
  author={Woronowicz, Stanis{\l }aw Lech},
  title={Compact quantum groups},
  conference={ title={Sym\'etries quantiques}, address={Les Houches}, date={1995}, },
  book={ publisher={North-Holland}, place={Amsterdam}, },
  date={1998},
  pages={845--884},
  review={\MRref {1616348}{99m:46164}},
}


\end{biblist}
\end{bibdiv}
\end{document}